\documentclass[13pt]{article}
\usepackage{latexsym}
\usepackage{geometry}
\usepackage{graphicx}
\usepackage{amsmath, amssymb, amsthm}
\usepackage{booktabs}
\usepackage{algorithm}
\usepackage{algorithmic}

\usepackage{url}
\usepackage{natbib}
\usepackage{appendix}
\usepackage{amsthm}
\usepackage{amsmath}

\usepackage{amsfonts}
\usepackage{multirow}
\usepackage{multicol}

\usepackage{color}
\usepackage{algorithm}
\usepackage{algorithmic}

\usepackage{xcolor,colortbl}
\definecolor{LightCyan}{rgb}{0.88,1,1}

\usepackage{graphicx}
\usepackage{subfigure}
\usepackage{hyperref}
\usepackage{tcolorbox}

\newtheorem{theorem}{Theorem}
\newtheorem{lemma}{Lemma}

\newtheorem{assumption}{Assumption}
\newtheorem{remark}{Remark}

\begin{document}
\title{ Near-Optimal Decentralized Momentum Method for Nonconvex-PL Minimax Problems }
\author{Feihu Huang\thanks{Feihu Huang is with College of Computer Science and Technology,
Nanjing University of Aeronautics and Astronautics, Nanjing, China;
and also with MIIT Key Laboratory of Pattern Analysis and Machine Intelligence, Nanjing, China. Email: huangfeihu2018@gmail.com}, \ Songcan Chen\thanks{Songcan Chen is with College of Computer Science and Technology, Nanjing University of Aeronautics and Astronautics, Nanjing, China; and also with MIIT Key Laboratory of Pattern Analysis and Machine Intelligence, Nanjing, China. E-mail: s.chen@nuaa.edu.cn} }

\date{}
\maketitle

\begin{abstract}
Minimax optimization plays an important role in many machine learning tasks such as generative adversarial
networks (GANs) and adversarial training. Although recently a wide variety of optimization methods have been proposed to solve the minimax problems, most of them
ignore the distributed setting where the data is distributed on multiple workers. Meanwhile, the existing decentralized minimax optimization methods rely on the strictly assumptions such as (strongly) concavity and
 variational inequality conditions. In the paper, thus,
we propose an efficient decentralized momentum-based gradient descent ascent (DM-GDA) method for the distributed nonconvex-PL minimax optimization, which is nonconvex in primal variable and is nonconcave in dual variable and satisfies the Polyak-Lojasiewicz (PL) condition. In particular, our DM-GDA method simultaneously uses the momentum-based techniques to update variables and estimate the stochastic gradients. Moreover, we provide a solid convergence analysis for our DM-GDA method, and prove that it obtains a near-optimal gradient complexity of $O(\epsilon^{-3})$ for finding an $\epsilon$-stationary solution of the nonconvex-PL stochastic minimax problems, which reaches the lower bound of nonconvex stochastic optimization. To the best of our knowledge, we first study the decentralized algorithm for Nonconvex-PL stochastic minimax optimization over a network.
\end{abstract}

\section{Introduction}
Minimax optimization has been widely applied in many machine learning tasks such as generative adversarial
networks (GANs)~\citep{goodfellow2014generative}, adversarial training~\citep{madry2018towards}, distributionally robust optimization~\citep{deng2020distributionally} and reinforcement
learning~\citep{zhang2021taming}. Meanwhile, decentralized distributed optimization~\citep{lian2017can} has received increasing attention in recent years in machine learning due to liberating the centralized agent with large
communication load and privacy risk.
In the paper, we consider the following distributed stochastic nonconvex-PL minimax problem over a communication network $G=(V,E)$, defined as
\begin{align} \label{eq:1}
 \min_{x \in \mathbb{R}^d} \max_{y \in \mathbb{R}^p} & \ f(x,y)\equiv\frac{1}{m}\sum_{i=1}^m \mathbb{E}_{\xi^i\sim \mathcal{D}^i}[f^i(x,y;\xi^i)],
\end{align}
where $f(x,y)\equiv\frac{1}{m}\sum_{i=1}^m f^i(x,y)$ and $f^i(x,y)=\mathbb{E}[f^i(x,y;\xi^i)]$
denotes the objective function in $i$-th client for any $i\in[m]$, which is differentiable and possibly nonconvex in primal variable $x$ and is differentiable and possibly nonconcave in dual variable $y$ and satisfies the Polyak-Lojasiewicz (PL) condition~\citep{polyak1963gradient,karimi2016linear}, which relaxes the strongly-convex in minimization optimization (i.e., strong-concave in maximization optimization). Here $\xi^i$ for any $i\in[m]$ is independent random variables follow unknown distributions $\mathcal{D}^i$, and for any $i,j\in [m]$ possibly $\mathcal{D}^i \neq \mathcal{D}^j$. $G=(V,E)$ is a communication network including $m$ computing clients, where
any agents $i,j\in V$ can communicate only if $(i,j)\in E$.

\begin{table*}
  \centering
  \caption{ \textbf{Gradient} (Stochastic First-order Oracle, i.e., SFO) complexity comparison of the representative \textbf{decentralized minimax optimization} algorithms for finding
  an $\epsilon$-stationary solution of the \textbf{nonconvex} minimax problems, i.e., $\mathbb{E}\|\nabla F(x)\|\leq \epsilon$
  or its equivalent variants, where $F(x)=\max_y \{f(x,y)\equiv\frac{1}{m}\sum_{i=1}^mf^i(x,y)\}$. Here $f(x,\cdot)$ denotes function \emph{w.r.t.} the second variable $y$ fixed $x$; $f(\cdot,y)$ denotes function \emph{w.r.t.} the first variable $x$ fixed $y$. \textbf{VI} denotes the variational inequality condition. $n$ denotes the sample size of \textbf{finite-sum} minimax problems, i.e., $\min_x \max_y \ \frac{1}{m} \sum_{i=1}^m \big( \frac{1}{n} \sum_{j=1}^n f^i_j(x,y)\big)$. Note that in the paper, we focus on the \textbf{stochastic} minimax optimization, as a finite-sum minimax optimization with the sample size $n\rightarrow + \infty$. }
  \label{tab:1}
   \resizebox{\textwidth}{!}{
\begin{tabular}{c|c|c|c|c}
  \hline
   \textbf{Algorithm} & \textbf{Reference} & \textbf{Gradient Complexity} & \textbf{Assumption } $f(\cdot,y)$& \textbf{Assumption } $f(x,\cdot)$  \\ \hline
  DPOSG  & \cite{liu2020decentralized}  & $O(\epsilon^{-12})$ & VI & VI \\  \hline
  GT/DA  & \cite{tsaknakis2020decentralized} & $\tilde{O}(n\epsilon^{-2})$ & Nonconvex & Concave \\  \hline
  DM-HSGD  & \cite{xian2021faster}  & $O(\epsilon^{-3})$ & Nonconvex & Strongly-Concave \\  \hline
  GT-SRVRI  & \cite{zhang2021taming}  & $O(\sqrt{n}\epsilon^{-2})$ & Nonconvex & Strongly-Concave \\  \hline
  DSGDA  & \cite{gao2022decentralized}  & $O(\sqrt{n}\epsilon^{-2})$ & Nonconvex & Strongly-Concave \\  \hline
  DREAM  & \cite{chen2022simple}  & $O(n+\sqrt{n}\epsilon^{-2})$/ $O(\epsilon^{-3})$ & Nonconvex & Strongly-Concave \\  \hline
  DM-GDA & Ours & $O(\epsilon^{-3})$ & Nonconvex & PL \\  \hline
\end{tabular}
 }
\end{table*}

When $m=1$, the problem~(\ref{eq:1}) reduces to a standard (stochastic) minimax optimization. Recently, many optimization methods~\citep{lin2019gradient,luo2020stochastic,huang2022accelerated} have been developed to solve these (stochastic) minimax problems. For example, \cite{lin2019gradient} studied the two-timescale (stochastic) Gradient Descent Ascent (GDA) methods for nonconvex (strongly) concave minimax optimization. Subsequently, \cite{luo2020stochastic,huang2022accelerated,huang2023adagda} proposed the efficient accelerated stochastic GDA methods for stochastic NonConvex-Strongly-Concave (NC-SC) minimax optimization. Meanwhile,
\cite{lu2020hybrid,chen2021proximal,huang2021efficient} studied the two-timescale proximal GDA methods
for the NC-SC minimax optimization with nonsmooth regularization.
Recently, some works have begun to studying the complex NonConvex NonConcave (NC-NC) minimax optimization.
For instance, \cite{yang2020global} presented an alternating GDA algorithm to solve minimax problems satisfying the  two-sided PL condition (i.e., PL-PL), which can linearly
converge to the stationary points. \cite{nouiehed2019solving} studied an effective
multi-step GDA method to solve the Nonconvex-PL minimax problems with one-sided PL condition. Subsequently, \cite{yang2022faster} developed
some efficient smoothed alternating GDA methods to solve Nonconvex-PL minimax problems. Meanwhile,
\cite{chen2022faster} studied a class of accelerated stochastic GDA methods based on the variance reduced technique for nonconvex minimax optimization with PL condition. More recently, \cite{huang2023enhanced} proposed a faster (adaptive) stochastic GDA method based on the momentum and variance reduced techniques, which reaches the near-optimal gradient complexity without large batches.

When $m\geq 2$, the problem~(\ref{eq:1}) is a standard distributed minimax optimization. Recently, some optimization methods~\citep{deng2021local,huang2022adaptive,tsaknakis2020decentralized,zhang2021taming,xian2021faster} have been studied in solving these distributed minimax problems. For example, \cite{deng2021local,sharma2022federated,huang2022adaptive} studied the federated learning methods for minimax optimization. Meanwhile, decentralized optimization methods~\citep{koppel2015saddle,mateos2015distributed,liu2020decentralized,beznosikov2021near,
rogozin2021decentralized,tsaknakis2020decentralized,zhang2021taming,xian2021faster} have been developed to solve the distributed minimax problems.
Specifically, \cite{koppel2015saddle,mateos2015distributed} proposed the decentralized algorithms for
convex-concave minimax optimization.
Subsequently, \cite{tsaknakis2020decentralized} studied the decentralized optimization methods for the NC-SC minimax optimization. Subsequently, \cite{xian2021faster} proposed a
faster decentralized minimax optimization method for NC-SC minimax optimization. \cite{liu2020decentralized} presented a fast decentralized parallel SGD method for a class of NC-NC minimax problems with the variational inequality conditions.
However, these decentralized minimax optimization methods rely on the strictly assumptions such as (strongly) concavity or variational inequality~\cite{iusem2017extragradient,liu2020decentralized} conditions. Clearly, these minimax methods can not well be competent to the general nonconvex-nonconcave minimax problems such as adversarial training deep neural networks (DNNs)~\citep{madry2018towards}.

In the paper, thus, we propose an efficient decentralized momentum-based gradient descent ascent (DM-GDA) method to the distributed nonconvex-PL minimax problem~(\ref{eq:1}), which is nonconvex in primal variable and is nonconcave in dual varaible and satisfies the PL condition.
Specifically, our DM-GDA algorithm simultaneously uses the momentum-based techniques to update variables and estimate the stochastic gradients.
In summary, our main contributions are:
\begin{itemize}
\item[(1)] We develop an efficient decentralized momentum-based gradient descent ascent (DM-GDA) method to the distributed nonconvex-PL minimax problem~(\ref{eq:1})
based on the momentum and variance reduced techniques. In particular, our DM-GDA algorithm simultaneously uses the momentum-based techniques to update variables and estimate the stochastic gradients.
\item[(2)] We introduce a solid convergence analysis framework for our DM-GDA algorithm, and prove that it reaches a near-optimal gradient (i.e., SFO) complexity of $O(\epsilon^{-3})$
for finding an $\epsilon$-stationary solution of Problem~(\ref{eq:1}).
\item[(3)] Numerical experimental results verify the efficiency of our DM-GDA algorithm on some robust DNNs training.
\end{itemize}
Note that since \cite{arjevani2022lower} provides a lower bound of gradient complexity $O(\epsilon^{-3})$ for finding an $\epsilon$-stationary solution of the nonconvex stochastic optimization $\min_x F(x)\equiv\mathbb{E}_{\xi}[f(x;\xi)]$, our gradient complexity of $O(\epsilon^{-3})$ is near optimal in solving
nonconvex-PL stochastic minimax optimization problem~(\ref{eq:1}). 

\section{Related Works}
In this section, we review some decentralized mini and minimax optimization algorithms, respectively.

\subsection{ Decentralized Mini Optimization }
Decentralized optimization has been widely applied to train large-scale machine learning tasks due to its efficiency and robustness in communication, which is no central
server with a communication bottleneck. Thus, recently numerous decentralized optimization methods~\citep{lian2017can,tang2018d,lu2019gnsd,koloskova2020unified,
pan2020d} have been proposed. For example, \cite{lian2017can} proposed an efficient decentralized SGD (i.e., D-SGD) method and theoretically proved the potential advantage of decentralized algorithm.
Subsequent, \cite{tang2018d} proposed an improved D-SGD method by relaxing a dissimilarity assumption.
\cite{lu2019gnsd,pu2021distributed} proposed efficient
 gradient-tracking-based D-SGD methods and investigated its convergence rates.
\cite{koloskova2020unified} studied the convergence
properties of decentralized SGD based on local updating with changing topology.

To further improve the gradient complexity of the decentralized SGD methods, a line of
workers~\citep{li2020communication,sun2020improving,xin2020variance,zhang2021gt,xin2021hybrid} studied the accelerated decentralized optimization methods based on the variance-reduced techniques.
For instance, \cite{xin2022fast} proposed some accelerated decentralized SGD methods based on
the variance-reduced techniques of SAGA~\citep{defazio2014saga} and SVRG~\citep{johnson2013accelerating}.
\cite{sun2020improving,pan2020d} presented the accelerated decentralized gradient-tracking methods by using
the variance reduced technique of SPIDER~\citep{fang2018spider}.
\cite{zhang2021gt,xin2021hybrid} developed the
decentralized variance-reduced SGD methods for nonconvex stochastic optimization based
on the gradient estimator of STORM~\citep{cutkosky2019momentum} and ProxHSGD~\citep{tran2022hybrid}.

To further improve the communication complexity of the decentralized
optimization methods, some communication-efficiency
decentralized methods~\citep{tang2018communication,koloskova2019decentralized,tang2019deep,reisizadeh2019robust} have been developed based on compressed and quantized
techniques. For example, \cite{koloskova2019decentralized} proposed an communication-efficiency gossip-based SGD algorithm by compressing (e.g., quantize or sparsify) the model
updates. \cite{reisizadeh2019robust} presented an communication-efficiency
decentralized gradient-based algorithm based on quantized message-passing.

\subsection{ Decentralized Minimax Optimization }
Minimax optimization is widely used in many machine learning tasks such as robust learning, AUC maximization
and reinforcement learning. In the big data settings, recently, there exists an increasing interest in distributed minimax optimization, e.g., adversarial training DNNs under distributed settings~\citep{liu2020decentralized}
 and policy evaluation over multi-agents~\citep{zhang2021taming}.
Recently, thus, decentralized optimization methods~\citep{koppel2015saddle,mateos2015distributed,liu2020decentralized,beznosikov2021near,
rogozin2021decentralized,tsaknakis2020decentralized,zhang2021taming,xian2021faster,gao2022decentralized,chen2022simple} have been developed to solve these distributed minimax optimization.
For example, \cite{koppel2015saddle,mateos2015distributed} developed the decentralized algorithms for
convex-concave minimax optimization.
Subsequently, \cite{tsaknakis2020decentralized} studied the decentralized optimization methods for the nonconvex-concave minimax optimization. Subsequently, \cite{xian2021faster} proposed a
faster decentralized minimax optimization method for nonconvex-strongly-concave stochastic minimax optimization.

\section{Preliminaries}

\subsection{Notations}
$[m]$ denotes the set $\{1,2,\cdots,m\}$.
$\|\cdot\|$ denotes the $\ell_2$ norm for vectors and spectral norm for matrices.
$\langle x,y\rangle$ denotes the inner product of two vectors $x$ and $y$. $I_{d}$ denotes a $d$-dimensional identity matrix. \textbf{1} is a vector of all one. 
Given function $f(x,y)$, $f(x,\cdot)$ denotes  function \emph{w.r.t.} the second variable with fixing $x$,
and $f(\cdot,y)$ denotes function \emph{w.r.t.} the first variable
with fixing $y$.
$a_m=O(b_m)$ denotes that $a_m \leq c b_m$ for some constant $c>0$. The notation $\tilde{O}(\cdot)$ hides logarithmic terms. Let $\mathbb{E}_t = \mathbb{E}_{\xi_t,\xi_{t-1},\cdots,\xi_1}$, where $\xi_t=\{\xi_t^i\}_{i=1}^m$. 

Let $\bar{x}_t=\frac{1}{m}\sum_{i=1}^m x^i_t$, $\bar{y}_t=\frac{1}{m}\sum_{i=1}^m y^i_t$, $\overline{\nabla_x f(x_t,y_t;\xi)}=\frac{1}{m}\sum_{i=1}^m \nabla_x f^i(x_t,y_t;\xi^i)$, $\overline{\nabla_x f(x_t,y_t)}=\frac{1}{m}\sum_{i=1}^m \nabla_x f^i(x_t,y_t)$, $\overline{\nabla_y f(x_t,y_t;\xi)}=\frac{1}{m}\sum_{i=1}^m \nabla_y f^i(x_t,y_t;\xi^i)$ and $\overline{\nabla_y f(x_t,y_t)}=\frac{1}{m}\sum_{i=1}^m \nabla_y f^i(x_t,y_t)$.

\subsection{ Mild Assumptions}
In this subsection, we give some mild assumptions on the problem \eqref{eq:1}.

\begin{assumption} \label{ass:1}
The graph $G=(V,E)$ is connected and undirected, which
can be represented by a mixing matrix $W\in \mathbb{R}^{m\times m}$: 1) $W_{i,j}>0$ if
$W_{i,j}\in E$ and $W_{i,j}=0$ otherwise; 2)
$W$ is doubly stochastic such that $W=W^T$, $\sum_{i=1}^mW_{i,j}=1$ and $\sum_{j=1}^mW_{i,j}=1$;
3) Eigenvalues of $W$ satisfy $\lambda_m \leq \cdots \leq \lambda_2 < \lambda_1=1$ and
$\nu=\max(|\lambda_2|, |\lambda_m|)<1$.
\end{assumption}

\begin{assumption} \label{ass:2}
For any $i \in [m]$, each component stochastic function $f^i(x,y;\xi^i)$ is $L_f$-smooth, such that
for all $x,x_1,x_2 \in \mathbb{R}^d$ and $y,y_1,y_2 \in \mathbb{R}^p$
\begin{align}
   & \|\nabla_x f^i(x_1,y;\xi^i)-\nabla_x f^i(x_2,y;\xi^i)\| \leq L_f\|x_1-x_2\|, \ \|\nabla_x f^i(x,y_1;\xi^i)-\nabla_x f^i(x,y_2;\xi^i)\| \leq L_f\|y_1-y_2\| \nonumber \\
   & \|\nabla_y f^i(x_1,y;\xi^i)-\nabla_y f^i(x_2,y;\xi^i)\| \leq L_f\|x_1-x_2\|, \ \|\nabla_y f^i(x,y_1;\xi^i)-\nabla_y f^i(x,y_2;\xi^i)\| \leq L_f\|y_1-y_2\|. \nonumber
\end{align}
\end{assumption}

\begin{assumption} \label{ass:3}
Stochastic function $f^i(x,y;\xi^i)$ has unbiased component stochastic gradient
with bounded variance for any $i\in [m]$, i.e.,
\begin{align}
& \mathbb{E}[\nabla_x f^i(x,y;\xi^i)] = \nabla_x f^i(x,y), \ \mathbb{E}\|\nabla_x f^i(x,y;\xi^i) - \nabla_x f^i(x,y) \|^2 \leq \sigma^2, \nonumber \\
& \mathbb{E}[\nabla_y f^i(x,y;\xi^i)] = \nabla_y f^i(x,y), \ \mathbb{E}\|\nabla_y f^i(x,y;\xi^i) - \nabla_y f^i(x,y) \|^2 \leq \sigma^2.
\end{align}
\end{assumption}

\begin{assumption} \label{ass:4}
The objective function $F(x)=\max_{y\in \mathbb{R}^p} \frac{1}{m}\sum_{i=1}^mf^i(x,y)$ is lower bounded, i.e.,  $F^* = \inf_{x\in
\mathbb{R}^d} F(x)$.
\end{assumption}

Assumption~\ref{ass:1} shows the protocol properties of connected network $G=(V,E)$, which is commonly used
in the decentralized optimization~\citep{lian2017can,koloskova2020unified}. Specifically, the mixing matrix $W$ represents the weights of averaging among the communication network topology, 
which is doubly stochastic, i.e., $W\textbf{1}=W^T\textbf{1}=\textbf{1}$, and 
its spectral gap satisfies $\|W-\textbf{1}\textbf{1}^T/m\|_2 = \nu \in (0,1)$.
Assumption~\ref{ass:3} shows the unbiased stochastic gradients $\nabla_x f^i(x,y;\xi^i)$ and $\nabla_y f^i(x,y;\xi^i)$. Assumption~\ref{ass:4} ensures the feasibility of the Problem~(\ref{eq:1}).
Assumptions~\ref{ass:2}-\ref{ass:4} are commonly used in stochastic minimax optimization problems \citep{nouiehed2019solving,huang2023enhanced}. 

Assumption~\ref{ass:2} shows the smoothness of stochastic functions $f^i(x,y;\xi^i)$ for all $i\in [m]$.
Based on Assumptions~\ref{ass:2} and~\ref{ass:3}, we have
\begin{align}
\|\nabla_x f^i(x_1,y)-\nabla_x f^i(x_2,y)\| & = \|\mathbb{E}[\nabla_x f^i(x_1,y;\xi^i)]-\mathbb{E}[\nabla_x f^i(x_2,y;\xi^i)]\| \nonumber \\
& \leq \mathbb{E}[\|\nabla_x f^i(x_1,y;\xi^i)-\nabla_x f^i(x_2,y;\xi^i)\| \leq L_f\|x_1-x_2\|,
\end{align}
and it is similar for $\|\nabla_x f^i(x,y_1)-\nabla_x f^i(x,y_2)\|\leq L_f\|y_1-y_2\|$,
$\|\nabla_y f^i(x,y_1)-\nabla_y f^i(x,y_2)\|\leq L_f\|y_1-y_2\|$ and $\|\nabla_y f^i(x_1,y)-\nabla_y f^i(x_2,y)\|\leq L_f\|x_1-x_2\|$. In other words, the function $f^i(x,y)$ for all $i\in [m]$ is $L_f$ smooth.
Similarly, let $f(x,y)=\frac{1}{m}\sum_{i=1}^mf^i(x,y)$, we have
\begin{align}
\|\nabla_x f(x_1,y) - \nabla_x f(x_2,y)\| & = \|\frac{1}{m}\sum_{i=1}^m\nabla_x f^i(x_1,y) - \frac{1}{m}\sum_{i=1}^m\nabla_x f^i(x_2,y)\| \nonumber \\
& \leq \frac{1}{m}\sum_{i=1}^m\|\nabla_x f^i(x_1,y) - \nabla_x f^i(x_2,y)\| \leq L_f\|x_1-x_2\|,
\end{align}
it is similar for $\|\nabla_x f(x,y_1)-\nabla_x f(x,y_2)\|\leq L_f\|y_1-y_2\|$,
$\|\nabla_y f(x,y_1)-\nabla_y f(x,y_2)\|\leq L_f\|y_1-y_2\|$ and $\|\nabla_y f(x_1,y)-\nabla_y f(x_2,y)\|\leq L_f\|x_1-x_2\|$. In other words, the function $f(x,y)$ is $L_f$ smooth.

\subsection{ Useful Lemmas }
In this subsection, we give some useful lemmas based on the above assumptions.

We first let $F(x)=\frac{1}{m}\sum_{i=1}^mF^i(x)$, $f(x,y) = \frac{1}{m}\sum_{i=1}^m f^i(x,y)$ and $F^i(x)=f^i(x,y^*(x))$ with $y^*(x)=\max_y \frac{1}{m}\sum_{i=1}^m f^i(x,y)$.

\begin{lemma} \label{lem:1}
Let $F(x)= f(x,y^*(x))$ with $y^*(x) \in \
\arg\max_y f(x,y)$. Under the above Assumptions~\ref{ass:1}-\ref{ass:2},
$\nabla F(x)=\nabla_x f(x,y^*(x))$ and $F(x)$ is $L$-smooth, i.e.,
\begin{align}
\|\nabla F(x_1) - \nabla F(x_2)\| \leq L\|x_1-x_2\|, \quad \forall x_1,x_2
\end{align}
where $L=L_f(1+\frac{\kappa}{2})$ with $\kappa=\frac{L_f}{\mu}$.
\end{lemma}

\begin{proof}
Based on Assumption \ref{ass:2}, we have
\begin{align}
& \|\nabla_x f^i(x_1,y)-\nabla_x f^i(x_2,y)\|=\|\mathbb{E}[\nabla_x f^i(x_1,y;\xi^i)-\nabla_x f^i(x_2,y;\xi^i)]\| \nonumber \\
& \leq \mathbb{E}\|\nabla_x f^i(x_1,y;\xi^i)-\nabla_x f^i(x_2,y;\xi^i)\| \leq L_f\|x_1-x_2\|,
\end{align}
then we can obtain $\nabla_x f^i(x,y)$ is $L_f$-Lipschitz continuous.
It is similar for $\nabla_y f^i(x,y)$.

According to Lemma A.5 of \cite{nouiehed2019solving}, then we have
\begin{align}
 \|\nabla F^i(x_1) - \nabla F^i(x_2)\| \leq L\|x_1-x_2\|, \quad \forall x_1,x_2
\end{align}
where $L=L_f(1+\frac{\kappa}{2})$ with $\kappa=\frac{L_f}{\mu}$.
Then we have
\begin{align}
 \|\nabla F(x_1)- \nabla F(x_2)\| & = \|\frac{1}{m}\sum_{i=1}^m\big(\nabla f^i(x_1)-\nabla f^i(x_2)\big)\| \nonumber \\
 &\leq
 \frac{1}{m}\sum_{i=1}^m\|\nabla f^i(x_1)-\nabla f^i(x_2)\| \leq L\|x_1-x_2\|,
\end{align}
i.e., the global function $F(x)$ also is $L$-smooth.
\end{proof}

\section{ Decentralized Nonconvex-PL Minimax Optimization }
In this section, we propose an novel efficient decentralized momentum-based gradient descent ascent (DM-GDA) to
solve Problem \eqref{eq:1}. Algorithm~\ref{alg:1} shows the detailed procedure of our DM-GDA algorithm.

\begin{algorithm}[t]
\caption{ DM-GDA Algorithm for Nonconvex-PL Minimax Optimization}
\label{alg:1}
\begin{algorithmic}[1]
\STATE {\bfseries Input:} $T>0$, tuning parameters $\{\gamma, \lambda, \eta_t, \alpha_t, \beta_t \}$; \\
\STATE {\bfseries initialize:} Set $x^i_0=\tilde{x}^i_0$ and $y^i_0=\tilde{y}^i_0$ for $i \in [m]$, and draw $m$ independent samples $\{\xi_{0}^i\}_{i=1}^m$
and then compute $u^i_{x,0} = \nabla_x f^i(x^i_0, y^i_0;\xi^i_{0})$, $u^i_{y,0} = \nabla_y f^i(x^i_0,y^i_0;\xi^i_{0})$, $w^i_{x,0} = \sum_{j\in \mathcal{N}_i} W_{i,j}u^j_{x,0}$ and $w^i_{y,0} = \sum_{j\in \mathcal{N}_i} W_{i,j}u^j_{y,0}$ for all $i \in [m]$. \\
\FOR{$t=0$ \textbf{to} $T-1$}
\FOR{$i=1,\cdots,m$ (\textbf{in parallel}) }
\STATE $\tilde{x}^i_{t+1} = \sum_{j\in \mathcal{N}_i} W_{i,j}x^j_t - \gamma w^i_{x,t}$ and $\tilde{y}^i_{t+1} = \sum_{j\in \mathcal{N}_i} W_{i,j}y^j_t + \lambda w^i_{y,t}$; \\
\STATE $x^i_{t+1} =  x^i_t + \eta_t(\tilde{x}^i_{t+1}-x^i_t) $ and $y^i_{t+1} =  y^i_t + \eta_t(\tilde{y}^i_{t+1}-y^i_t) $ ;  \\
\STATE Draw a sample $\xi^i_{t+1}\sim \mathcal{D}^i$;
\STATE $u^i_{x,t+1} = \nabla_x f^i(x^i_{t+1},y^i_{t+1};\xi^i_{t+1}) + (1-\alpha_{t+1})(u^i_{x,t} -\nabla_x f^i(x^i_t,y^i_t;\xi^i_{t+1}))$; \\
\STATE $u^i_{y,t+1} = \nabla_y f^i(x^i_{t+1},y^i_{t+1};\xi^i_{t+1}) + (1-\beta_{t+1})(u^i_{y,t} -\nabla_y f^i(x^i_t,y^i_t;\xi^i_{t+1}))$; \\
\STATE $w^i_{x,t+1} =  \sum_{j\in \mathcal{N}_i}W_{i,j}\big(w^j_{x,t} + u^j_{x,t+1}-u^j_{x,t} \big)$; \\
\STATE $w^i_{y,t+1} =  \sum_{j\in \mathcal{N}_i}W_{i,j}\big(w^j_{y,t} + u^j_{y,t+1}-u^j_{y,t} \big)$; \\
\ENDFOR
\ENDFOR
\STATE {\bfseries Output:} Chosen uniformly random from $\{x^i_{t\geq1},y^i_{t\geq1}\}_{i=1}^m$.
\end{algorithmic}
\end{algorithm}

At the line 5 of Algorithm~\ref{alg:1}, each client uses the stochastic gradient descent to update the local primal variable and
the stochastic gradient ascent to update the local dual variable:
\begin{align}
 \tilde{x}^i_{t+1} = \sum_{j\in \mathcal{N}_i} W_{i,j}x^j_t - \gamma w^i_{x,t}, \quad \tilde{y}^i_{t+1} = \sum_{j\in \mathcal{N}_i} W_{i,j}y^j_t + \lambda w^i_{y,t},
\end{align}
where the constants $\gamma$ and $\lambda$ are the step-sizes for individual primal and dual variables,  respectively. Here $\mathcal{N}_i=\{j\in V \ | \ (i,j)\in E, j=i\}$
denotes the neighborhood of the $i$-th client.  Here each client communicates with their neighbors
to update the variables $x$ and $y$. Then we use the momentum iteration to simultaneously update the variables $x$ and $y$ at the line 6 of Algorithm~\ref{alg:1}.
In lines 8-9 of Algorithm \ref{alg:1}, each client uses the momentum-based variance reduced technique of STORM~\citep{cutkosky2019momentum}/ ProxHSGD~\citep{tran2022hybrid}
to update the stochastic gradients based on local data: for $i\in [m]$
\begin{align}
& u^i_{x,t+1} = \nabla_x f^i(x^i_{t+1};\xi^i_{t+1}) + (1-\alpha_{t+1})(u^i_{x,t} -\nabla_x f^i(x^i_{t};\xi^i_{t+1}))  \\
& u^i_{y,t+1} = \nabla_y f^i(x^i_{t+1};\xi^i_{t+1}) + (1-\beta_{t+1})(u^i_{y,t} -\nabla_y f^i(x^i_{t};\xi^i_{t+1})),
\end{align}
where $\alpha_{t+1}, \beta_{t+1}\in (0,1)$.
At the lines 10-11 of our Algorithm \ref{alg:1}, then each client communicates with its neighbors
to compute gradient estimators $w^i_{x,t+1}$ and $w^i_{y,t+1}$, defined as 
\begin{align}
& w^i_{x,t+1} =  \sum_{j\in \mathcal{N}_i}W_{i,j}\big(w^j_{x,t} + u^j_{x,t+1}-u^j_{x,t} \big),  \\
& w^i_{y,t+1} =  \sum_{j\in \mathcal{N}_i}W_{i,j}\big(w^j_{y,t} + u^j_{y,t+1}-u^j_{y,t} \big).
\end{align}
Here we use gradient tracking technique~\citep{xu2015augmented,di2016next} to reduce the consensus error. 
Specifically, local gradients $w^i_{x,t+1}$ and $w^i_{y,t+1} $ track the directions of global gradients.

\section{Convergence Analysis}
In this section, we provide the convergence properties of our DM-GDA algorithm
under some mild assumptions. All related proofs are provided in the following Appendix.
We first review some useful lemmas and assumptions.

\begin{lemma} \label{lem:2}
 The sequences $\{u^i_{x,t},u^i_{y,t},w^i_{x,t},w^i_{y,t}\}_{i=1}^m$ be generated from our Algorithm~\ref{alg:1}, we have for all $t\geq 1$,
 \begin{align}
 \frac{1}{m}\sum_{i=1}^mu^i_{x,t} = \bar{u}_{x,t}= \bar{w}_{x,t} =\frac{1}{m}\sum_{i=1}^mw^i_{x,t}, \quad \frac{1}{m}\sum_{i=1}^mu^i_{y,t} = \bar{u}_{y,t}= \bar{w}_{y,t} =\frac{1}{m}\sum_{i=1}^mw^i_{y,t}.
 \end{align}
\end{lemma}

\begin{lemma} \label{lem:3}
Suppose the sequence $\{\bar{x}_t,\bar{y}_t\}_{t=1}^T$ be generated from Algorithm~\ref{alg:1}, where 
$\bar{x}_t = \frac{1}{m}\sum_{i=1}^mx^i_t$ and $\bar{y}_t = \frac{1}{m}\sum_{i=1}^my^i_t$. 
Under the Assumptions~\ref{ass:1}-\ref{ass:2}, given $0<\gamma\leq \frac{\lambda\mu}{16L}$ and $0<\lambda \leq \frac{1}{2L_f\eta_t}$ for all $t\geq 1$, we have
\begin{align}
F(\bar{x}_{t+1}) - f(\bar{x}_{t+1},\bar{y}_{t+1})
& \leq (1-\frac{\eta_t\lambda\mu}{2}) \big(F(\bar{x}_t) -f(\bar{x}_t,\bar{y}_t)\big) + \frac{\eta_t\gamma}{8}\|\bar{w}_{x,t}\|^2  -\frac{\eta_t\lambda}{4}\|\bar{w}_{y,t}\|^2 \nonumber \\
& \quad + \eta_t\lambda\|\nabla_y f(\bar{x}_t,\bar{y}_t)-\bar{w}_{y,t}\|^2,
\end{align}
where $F(\bar{x}_t)=f(\bar{x}_t,y^*(\bar{x}_t))$ with $y^*(\bar{x}_t) \in \arg\max_{y}f(\bar{x}_t,y)$ for all $t\geq 1$.
\end{lemma}

Lemma~\ref{lem:3} shows property of the residual $F(x)-f(x,y) = \max_yf(x,y) -f(x,y) \geq 0$, as Lemma 1 of \cite{huang2023enhanced}. 
We define a useful Lyapunov function (i.e., potential function), for any $t\geq 1$
\begin{align} 
\Omega_t & = \mathbb{E}_t \Big[ F(\bar{x}_t) + \frac{72\gamma L^2_f}{\lambda\mu^2}\big(F(\bar{x}_t) -f(\bar{x}_t,\bar{y}_t)\big) + (\rho_{x,t-1}-\frac{9\gamma\eta}{2})\|\bar{u}_{x,t-1} - \overline{\nabla_x f(x_{t-1},y_{t-1})}\|^2  \nonumber \\
& \quad + (\rho_{y,t-1}-\frac{144\gamma\eta L^2_f}{\mu^2})\|\bar{u}_{y,t-1} - \overline{\nabla_y f(x_{t-1},y_{t-1})}\|^2 + \varrho_{x,t-1}\frac{1}{m}\sum_{i=1}^m \|u^i_{x,t-1} - \nabla_x f^i(x^i_{t-1},y^i_{t-1})\|^2   \nonumber \\
& \quad + \varrho_{y,t-1}\frac{1}{m}\sum_{i=1}^m \|u^i_{y,t-1} - \nabla_y f^i(x^i_{t-1},y^i_{t-1})\|^2 + (\theta_{x,t-1}-9\gamma\eta L_f^2-\frac{\gamma\eta L^2}{2}-\frac{288\gamma\eta L^4_f}{\mu^2})\frac{1}{m}\sum_{i=1}^m\|x^i_{t-1}-\bar{x}_{t-1}\|^2  \nonumber \\
& \quad  + (\theta_{y,t-1}-18\gamma\eta L_f^2-\frac{288\gamma\eta L^4_f}{\mu^2})\frac{1}{m}\sum_{i=1}^m\|y^i_{t-1}-\bar{y}_{t-1}\|^2 +(\vartheta_{x,t-1}-\frac{3\gamma\eta}{4})\frac{1}{m}\sum_{i=1}^m\|w^i_{x,t-1}-\bar{w}_{x,t-1}\|^2
 \nonumber \\
& \quad +\vartheta_{y,t-1}\frac{1}{m}\sum_{i=1}^m\|w^i_{y,t-1}-\bar{w}_{y,t-1}\|^2 +
\frac{\gamma\eta}{12}\frac{1}{m}\sum_{i=1}^m\|w^i_{x,t-1}\|^2 + \frac{18\gamma L^2_f\eta}{\mu^2}\|\bar{w}_{y,t-1}\|^2\Big],
\end{align}
where $\rho_{x,t-1}\geq \frac{9\gamma\eta}{2}$, $\rho_{y,t-1}\geq \frac{144\gamma\eta L^2_f}{\mu^2}$, $\varrho_{x,t-1}> 0$, $\varrho_{y,t-1}> 0$,
$\theta_{x,t-1}\geq 9\gamma\eta L_f^2+\frac{\gamma\eta L^2}{2}+\frac{288\gamma\eta L^4_f}{\mu^2}$, $\theta_{y,t-1}\geq 18\gamma\eta L_f^2+\frac{288\gamma\eta L^4_f}{\mu^2}$,
$\vartheta_{x,t-1}\geq \frac{3\gamma\eta}{4}$ and $\vartheta_{y,t-1}> 0$ for all $t\geq1$.
Then based on the above potential function $\Omega_t$, we provide the convergence properties of 
our DM-GDA algorithm. 

\begin{theorem}  \label{th:1}
 Suppose the sequences $\{\bar{x}_t,\bar{y}_t\}_{t=1}^T$ be generated from Algorithm~\ref{alg:1}.
 Under the above Assumptions~\ref{ass:1}-\ref{ass:4}, and let $\alpha_t =\beta_t = O(\frac{1}{T^{2/3}})$ and $\eta_t=\eta=O(\frac{1}{T^{1/3}})$ for all $t\geq1$, $0<\gamma\leq \frac{\lambda\mu}{16L}$ and $\lambda =O(1)$ for all $t\geq 0$, we have
\begin{align}
  \frac{1}{T}\sum_{t=1}^T \mathbb{E}\|\nabla F(\bar{x}_t)\|
  \leq O(\frac{1}{T^{1/3}}+\frac{\sigma^2}{T^{1/3}}),
\end{align}
where $F(x)=\max_{y\mathbb{R}^p} f(x,y)=\max_{y\mathbb{R}^p}\frac{1}{m}\sum_{i=1}^m\mathbb{E}[f^i(x,y;\xi^i)]$. 
\end{theorem}

\begin{remark}
From the above Algorithm~\ref{alg:1}, our DM-GDA algorithm has a convergence rate $O(\frac{1}{T^{1/3}})$. 
Since our DM-GDA algorithm needs four stochastic gradients at each iteration, let $O(\frac{1}{T^{1/3}})\leq \epsilon$, 
it has a gradient complexity of $4 \cdot T = O(\epsilon^{-3})$.
Thus, our DM-GDA algorithm requires gradient complexity of $O(\epsilon^{-3})$ 
for finding an $\epsilon$-stationary solution of Problem (\ref{eq:1}). Since \cite{arjevani2022lower} has a lower bound of gradient complexity $O(\epsilon^{-3})$ for finding an $\epsilon$-stationary solution of the nonconvex stochastic optimization $\min_x F(x)\equiv\mathbb{E}_{\xi}[f(x;\xi)]$, our gradient complexity of $O(\epsilon^{-3})$ is near optimal. 
\end{remark}

\section{Numerical Experiments}
In this section, we use some numerical experiments...

\section{Conclusion}
In the paper, we studied a class of nonconvex nonconcave distributed minimax problems, where the dual variable satisfies the PL condition.
Then we proposed an efficient decentralized momentum-based gradient descent ascent (DM-GDA) method to solve these distributed nonconvex-PL minimax problems. Moreover, we provide a solid convergence analysis framework for our DM-GDA algorithm,
and proved that it reaches a near-optimal gradient complexity of $O(\epsilon^{-3})$ for finding an $\epsilon$-stationary solution of Problem~(\ref{eq:1}).


\small

\bibliographystyle{plainnat}

\bibliography{Minimax-NCPL}

\newpage

\appendix

\section{Appendix}
In this section, we provide the detailed convergence analysis of our DM-GDA algorithm.
We first review and provide some useful lemmas. 

For notational simplicity, let $x_t=[(x^1_t)^T,\cdots,(x^m_t)^T]^T\in \mathbb{R}^{md}$, $\tilde{x}_t=[(\tilde{x}^1_t)^T,\cdots,(\tilde{x}^m_t)^T]^T\in \mathbb{R}^{md}$, $y_t=[(y^1_t)^T,\cdots,(y^m_t)^T]^T\in \mathbb{R}^{mp}$, $\tilde{y}_t=[(\tilde{y}^1_t)^T,\cdots,(\tilde{y}^m_t)^T]^T\in \mathbb{R}^{mp}$, $w_{x,t}=[(w^1_{x,t})^T,\cdots,(w^m_{x,t})^T]^T\in \mathbb{R}^{md}$ and $w_{y,t}=[(w^1_{y,t})^T,\cdots,(w^m_{y,t})^T]^T\in \mathbb{R}^{mp}$ for all $t\geq1$.

\begin{lemma} \label{lem:A1}
(Lemma A.5 of \cite{nouiehed2019solving})
Let $F(x)= f(x,y^*(x))=\max_y f(x,y)$ with $y^*(x) \in \
\arg\max_y f(x,y)$. Under the above Assumptions~\ref{ass:1}-\ref{ass:2},
$\nabla F(x)=\nabla_x f(x,y^*(x))$ and $F(x)$ is $L$-smooth, i.e.,
\begin{align}
\|\nabla F(x_1) - \nabla F(x_2)\| \leq L\|x_1-x_2\|, \quad \forall x_1,x_2
\end{align}
where $L=L_f(1+\frac{\kappa}{2})$ with $\kappa=\frac{L_f}{\mu}$.
\end{lemma}

\begin{lemma} \label{lem:A2}
(\cite{karimi2016linear})
 Function $f(x): \mathbb{R}^d\rightarrow \mathbb{R}$ is $L$-smooth and satisfies PL condition with constant $\mu$, then it also
satisfies error bound (EB) condition with $\mu$, i.e., for all $x \in \mathbb{R}^d$
\begin{align}
 \|\nabla f(x)\| \geq \mu\|x^*-x\|,
\end{align}
where $x^* \in \arg\min_{x} f(x)$. It also satisfies quadratic growth (QG) condition with $\mu$, i.e.,
\begin{align}
 f(x)-\min_x f(x) \geq \frac{\mu}{2}\|x^*-x\|^2.
\end{align}
\end{lemma}
From the above lemma~\ref{lem:A2}, when consider the problem $\max_x f(x)$ that is equivalent to the problem $-\min_x -f(x)$, we have
\begin{align}
 & \|\nabla f(x)\| \geq \mu\|x^*-x\|, \\
 & \max_x f(x) - f(x) \geq \frac{\mu}{2}\|x^*-x\|^2.
\end{align}

\subsection{ Convergence Analysis of DM-GDA Algorithm }

\begin{lemma} \label{lem:B1}
(Restatement of Lemma 2)
 The sequences $\{u^i_{x,t},u^i_{y,t},w^i_{x,t},w^i_{y,t}\}_{i=1}^m$ be generated from our Algorithm~\ref{alg:1}, we have for all $t\geq 1$,
 \begin{align}
 \frac{1}{m}\sum_{i=1}^mu^i_{x,t} = \bar{u}_{x,t}= \bar{w}_{x,t} =\frac{1}{m}\sum_{i=1}^mw^i_{x,t}, \quad \frac{1}{m}\sum_{i=1}^mu^i_{y,t} = \bar{u}_{y,t}= \bar{w}_{y,t} =\frac{1}{m}\sum_{i=1}^mw^i_{y,t}.
 \end{align}
\end{lemma}

\begin{proof}
We proceed by induction.
From our Algorithm~\ref{alg:1}, since $w^i_{x,1} = \sum_{j\in \mathcal{N}_i} W_{i,j}u^j_{x,1}$, we have
\begin{align}
 \bar{w}_{x,1} = \frac{1}{m}\sum_{i=1}^m w^i_{x,1} = \frac{1}{m}\sum_{i=1}^m\sum_{j\in \mathcal{N}_i} W_{i,j}u^j_{x,1}
 =\frac{1}{m}\sum_{j=1}^m u^j_{x,1}\sum_{i=1}^m W_{i,j} = \frac{1}{m}\sum_{j=1}^m u^j_{x,1} = \bar{u}_{x,1},
\end{align}
where the second last equality is due to $\sum_{i=1}^m W_{i,j}=1$ from Assumption \ref{ass:1}.

From the line 10 of Algorithm~\ref{alg:1}, we have for all $t\geq 1$
\begin{align}
  \bar{w}_{x,t+1} & = \frac{1}{m}\sum_{i=1}^m w^i_{x,t+1} = \frac{1}{m}\sum_{i=1}^m\sum_{j\in \mathcal{N}_i} W_{i,j}\big(w^j_{x,t} + u^j_{x,t+1}-u^j_{x,t}\big) \nonumber \\
  & = \frac{1}{m}\sum_{j=1}^m\big(w^j_{x,t}
  +u^j_{x,t+1}-u^j_{x,t}\big)\sum_{i=1}^m W_{i,j} \nonumber \\
  & = \frac{1}{m}\sum_{j=1}^m\big(w^j_{x,t}
  +u^j_{x,t+1}-u^j_{x,t}\big) = \bar{w}_{x,t}
  +\bar{u}_{x,t+1}-\bar{u}_{x,t} = \bar{u}_{x,t+1},
\end{align}
where the second last equality is due to $\mathcal{N}_i=\{j\in V \ | \ (i,j)\in E, j=i\}$ and $\sum_{i=1}^m W_{i,j}=1$, and the last equality holds by
the inductive hypothesis, i.e., $\bar{w}_{x,t}=\bar{u}_{x,t}$.

Similarly, we can get
 \begin{align}
  \frac{1}{m}\sum_{i=1}^mu^i_{y,t} = \bar{u}_{y,t}= \bar{w}_{y,t} =\frac{1}{m}\sum_{i=1}^mw^i_{y,t}.
 \end{align}

\end{proof}

\begin{lemma} \label{lem:C1}
Under the above assumptions, and assume the stochastic gradient estimators $\big\{u^i_{x,t},u^i_{y,t}\big\}_{i=1}^m$ be generated from Algorithm~\ref{alg:1}, we have for $\forall i\in[m]$,
 \begin{align}
\mathbb{E}\|u^i_{x,t+1} - \nabla_x f^i(x^i_{t+1},y^i_{t+1})\|^2
  & \leq (1-\alpha_{t+1})\mathbb{E} \|u^i_{x,t} - \nabla_x f^i(x^i_t,y^i_t)\|^2 + 2\alpha_{t+1}^2\sigma^2 \nonumber \\
  & \quad + 2L^2_f\eta_t^2\mathbb{E}\big(\|\tilde{x}^i_{t+1}-x^i_t\|^2+\|\tilde{y}^i_{t+1}-y^i_t\|^2\big), \nonumber
 \end{align}
  \begin{align}
 \mathbb{E}\|\bar{u}_{x,t+1} - \overline{\nabla_x f(x_{t+1},y_{t+1})}\|^2
  & \leq (1-\alpha_{t+1})\mathbb{E} \|\bar{u}_{x,t} -\overline{\nabla_x f(x_t,y_t)}\|^2 + \frac{2\alpha_{t+1}^2\sigma^2}{m} \nonumber \\
  & \quad + \frac{2L^2_f\eta_t^2}{m^2}\sum_{i=1}^m\mathbb{E}\big(\|\tilde{x}^i_{t+1}-x^i_t\|^2+\|\tilde{y}^i_{t+1}-y^i_t\|^2\big), \nonumber
  \end{align}
 \begin{align}
 \mathbb{E}\|u^i_{y,t+1} - \nabla_y f^i(x^i_{t+1},y^i_{t+1})\|^2
  & \leq (1-\beta_{t+1})\mathbb{E} \|u^i_{y,t} - \nabla_y f^i(x^i_t,y^i_t)\|^2 + 2\beta_{t+1}^2\sigma^2 \nonumber \\
  & \quad + 2L^2_f\eta_t^2\mathbb{E}\big(\|\tilde{x}^i_{t+1}-x^i_t\|^2+\|\tilde{y}^i_{t+1}-y^i_t\|^2\big), \nonumber
 \end{align}
 \begin{align}
 \mathbb{E}\|\bar{u}_{y,t+1} - \overline{\nabla_y f(x_{t+1},y_{t+1})}\|^2
 & \leq (1-\beta_{t+1})\mathbb{E} \|\bar{u}_{y,t} -\overline{\nabla_y f(x_t,y_t)}\|^2 + \frac{2\beta_{t+1}^2\sigma^2}{m} \nonumber \\
 & \quad + \frac{2L^2_f\eta_t^2}{m^2}\sum_{i=1}^m\mathbb{E}\big(\|\tilde{x}^i_{t+1}-x^i_t\|^2+\|\tilde{y}^i_{t+1}-y^i_t\|^2\big), \nonumber
 \end{align}
where $\overline{\nabla_x f(x_t,y_t)} = \frac{1}{m}\sum_{i=1}^m\nabla_x f^i(x^i_t,y^i_t)$ and $\overline{\nabla_y f(x_t,y_t)} = \frac{1}{m}\sum_{i=1}^m\nabla_y f^i(x^i_t,y^i_t)$.
\end{lemma}

\begin{proof}
Since $u^i_{x,t+1} = \nabla_x f^i(x^i_{t+1},y^i_{t+1};\xi^i_{t+1}) + (1-\alpha_{t+1})(u^i_{x,t}-\nabla_x f^i(x^i_t,y^i_t;\xi^i_{t+1}))$ for any $i\in [m]$, we have
\begin{align}
\bar{u}_{x,t+1} & = \frac{1}{m}\sum_{i=1}^m\big( \nabla_x f^i(x^i_{t+1},y^i_{t+1};\xi^i_{t+1}) + (1-\alpha_{t+1})(u^i_{x,t}-\nabla_x f^i(x^i_t,y^i_t;\xi^i_{t+1}))\big)\nonumber \\
& =  \overline{\nabla_x f(x_{t+1},y_{t+1};\xi_{t+1})} + (1-\alpha_{t+1})(\bar{u}_{x,t} - \overline{\nabla_x f(x_t,y_t;\xi_{t+1})}).
\end{align}

Then we have
\begin{align}
 &\mathbb{E}\|\bar{u}_{x,t+1} - \overline{\nabla_x f(x_{t+1},y_{t+1})}\|^2  \\
 & = \mathbb{E}\|\frac{1}{m}\sum_{i=1}^m\big(u^i_{x,t+1} - \nabla_x f^i(x^i_{t+1},y^i_{t+1})\big)\|^2 \nonumber \\
 & = \mathbb{E}\|\frac{1}{m}\sum_{i=1}^m\big(\nabla_x f^i(x^i_{t+1},y^i_{t+1};\xi^i_{t+1}) + (1-\alpha_{t+1})\big(u^i_{x,t}
 - \nabla_x f^i(x^i_t,y^i_t;\xi^i_{t+1})\big) - \nabla_x f^i(x^i_{t+1},y^i_{t+1})\big)\|^2 \nonumber \\
 & = \mathbb{E}\big\|\frac{1}{m}\sum_{i=1}^m\Big( \nabla_x f^i(x^i_{t+1},y^i_{t+1};\xi^i_{t+1})-\nabla_x f^i(x^i_{t+1},y^i_{t+1}) - (1-\alpha_{t+1})\big( \nabla_x f^i(x^i_t,y^i_t;\xi^i_{t+1}) - \nabla_x f^i(x^i_t,y^i_t) \big) \Big) \nonumber \\
 & \quad + (1-\alpha_{t+1})\frac{1}{m}\sum_{i=1}^m\big( u^i_{x,t} - \nabla_x f^i(x^i_t,y^i_t) \big) \big\|^2 \nonumber \\
 & = \frac{1}{m^2}\sum_{i=1}^m\mathbb{E}\big\|\nabla_x f^i(x^i_{t+1},y^i_{t+1};\xi^i_{t+1})-\nabla_x f^i(x^i_{t+1},y^i_{t+1}) - (1-\alpha_{t+1})\big( \nabla_x f^i(x^i_t,y^i_t;\xi^i_{t+1}) - \nabla_x f^i(x^i_t,y^i_t) \big) \big\|^2 \nonumber \\
 & \quad + (1-\alpha_{t+1})^2\mathbb{E}\| \bar{u}_{x,t} - \overline{\nabla_x f(x_t,y_t)} \big\|^2 \nonumber \\
 & \leq \frac{2(1-\alpha_{t+1})^2}{m^2}\sum_{i=1}^m\mathbb{E}\|\nabla_x f^i(x^i_{t+1},y^i_{t+1};\xi^i_{t+1})-\nabla_x f^i(x^i_{t+1},y^i_{t+1}) - \nabla_x f^i(x^i_t,y^i_t;\xi^i_{t+1}) + \nabla_x f^i(x^i_t,y^i_t) \big) \|^2 \nonumber \\
 & \quad + \frac{2\alpha^2_{t+1}}{m^2}\sum_{i=1}^m\mathbb{E}\|\nabla_x f^i(x^i_{t+1},y^i_{t+1};\xi^i_{t+1})-\nabla_x f^i(x^i_{t+1},y^i_{t+1}) \|^2 + (1-\alpha_{t+1})^2\mathbb{E}\|\bar{u}_{x,t} - \overline{\nabla_x f(x_t,y_t)} \|^2 \nonumber \\
 & \leq \frac{2(1-\alpha_{t+1})^2}{m^2}\sum_{i=1}^m\mathbb{E}\|\nabla_x f^i(x^i_{t+1},y^i_{t+1};\xi^i_{t+1})- \nabla_x f^i(x^i_t,y^i_t;\xi^i_{t+1}) \|^2 + \frac{2\alpha^2_{t+1}\sigma^2}{m}
  + (1-\alpha_{t+1})^2\mathbb{E}\|\bar{u}_{x,t} - \overline{\nabla_x f(x_t,y_t)}\|^2 \nonumber \\
 & \leq (1-\alpha_{t+1})^2 \mathbb{E}\|\bar{u}_{x,t} - \overline{\nabla_x f(x_t,y_t)}\|^2 + \frac{2\alpha^2_{t+1}\sigma^2}{m} + \frac{2(1-\alpha_{t+1})^2L^2_f}{m^2}\sum_{i=1}^m\mathbb{E}
 \big(\|\tilde{x}^i_{t+1}-x^i_t\|^2+\|\tilde{y}^i_{t+1}-y^i_t\|^2\big)   \nonumber \\
 & \leq (1-\alpha_{t+1}) \mathbb{E}\|\bar{u}_{x,t} - \overline{\nabla_x f(x_t,y_t)} \|^2 + \frac{2\alpha^2_{t+1}\sigma^2}{m} + \frac{2L^2_f\eta^2_t}{m^2}\sum_{i=1}^m\mathbb{E}\big(\|\tilde{x}^i_{t+1}-x^i_t\|^2+\|\tilde{y}^i_{t+1}-y^i_t\|^2\big) \nonumber,
\end{align}
where the forth equality holds by the following fact: for any $i\in [m]$,
\begin{align}
 \mathbb{E}_{\xi^i_{t+1}} \big[\nabla_x f^i(x^i_{t+1},y^i_{t+1};\xi^i_{t+1})\big]=\nabla_x f^i(x^i_{t+1},y^i_{t+1}), \
 \mathbb{E}_{\xi^i_{t+1}} \big[\nabla_x f^i(x^i_t,y^i_t;\zeta^i_{t+1}))\big]= \nabla_x f^i(x^i_t,y^i_t), \nonumber
\end{align}
and for any $i \neq j\in [m]$, $\zeta^i_{t+1}$ and $\zeta^j_{t+1}$ are independent;
the second inequality holds by the inequality $\mathbb{E}\|\zeta-\mathbb{E}[\zeta]\|^2 \leq \mathbb{E}\|\zeta\|^2$ and Assumption \ref{ass:3};
the second last inequality is due to Assumption \ref{ass:2}; the last inequality holds by $0<\alpha_{t+1} \leq 1$ and $x^i_{t+1}=x^i_t+\eta_t(\tilde{x}^i_{t+1}-x^i_t)$.

Similarly, we have for $\forall i\in[m]$,
\begin{align}
 \mathbb{E}\|u^i_{x,t+1} - \nabla_x f^i(x^i_{t+1},y^i_{t+1})\|^2
  & \leq (1-\alpha_{t+1})\mathbb{E} \|u^i_{x,t} - \nabla_x f^i(x^i_t,y^i_t)\|^2 + 2\alpha_{t+1}^2\sigma^2 \nonumber \\
  & \quad + 2L^2_f\eta_t^2\mathbb{E}\big(\|\tilde{x}^i_{t+1}-x^i_t\|^2+\|\tilde{y}^i_{t+1}-y^i_t\|^2\big).
 \end{align}

\end{proof}

\begin{lemma} \label{lem:D1}
Given the sequences $\big\{x^i_{t},y^i_{t},w^i_{x,t},w^i_{y,t}\big\}_{i=1}^m$ be generated from Algorithm~\ref{alg:1}. We have
\begin{align}
\sum_{i=1}^m\|x^i_{t+1} - \bar{x}_{t+1}\|^2 & \leq (1-\frac{(1-\nu^2)\eta_t}{2})\sum_{i=1}^m\|x^i_t-\bar{x}_t\|^2
+ \frac{2\eta_t\gamma^2}{1-\nu^2}\sum_{i=1}^m\| w^i_{x,t}-\bar{w}_{x,t} \|^2,
\nonumber \\
\sum_{i=1}^m\|\tilde{x}^i_{t+1}-x^i_t\|^2
 & \leq (3+\nu^2)\sum_{i=1}^m\|x^i_t -\bar{x}_t\|^2 + \frac{2(1+\nu^2)}{1-\nu^2}\gamma^2
 \sum_{i=1}^m\|w^i_{x,t}\|^2,  \nonumber \\
\sum_{i=1}^m\|w^i_{x,t+1}-\bar{w}_{x,t+1}\|^2
    & \leq \nu\sum_{i=1}^m\|w^i_{x,t}-\bar{w}_{x,t}\|^2 + \frac{\nu^2}{1-\nu}\big(4\alpha_{t+1}^2\sum_{i=1}^m\|u^i_{x,t}-\nabla_x f^i(x^i_t,y^i_t)\|^2 \nonumber \\
    & \quad + 4\alpha^2_{t+1}m\sigma^2 + 16\eta^2_tL_f^2\sum_{i=1}^m\big( \|\tilde{x}^i_{t+1}-x^i_t\|^2 + \|\tilde{y}^i_{t+1}-y^i_t\|^2\big) \big), \nonumber \\
\sum_{i=1}^m\|y^i_{t+1} - \bar{y}_{t+1}\|^2 & \leq (1-\frac{(1-\nu^2)\eta_t}{2})\sum_{i=1}^m\|y^i_t-\bar{y}_t\|^2
+ \frac{2\eta_t\lambda^2}{1-\nu^2}\sum_{i=1}^m\| w^i_{y,t}-\bar{w}_{y,t} \|^2,
\nonumber \\
\sum_{i=1}^m\|\tilde{y}^i_{t+1}-y^i_t\|^2
 & \leq (3+\nu^2)\sum_{i=1}^m\|y^i_t -\bar{y}_t\|^2 + \frac{2(1+\nu^2)}{1-\nu^2}\lambda^2
 \sum_{i=1}^m\|w^i_{y,t}\|^2,  \nonumber \\
\sum_{i=1}^m\|w^i_{y,t+1}-\bar{w}_{y,t+1}\|^2
    & \leq \nu\sum_{i=1}^m\|w^i_{y,t}-\bar{w}_{y,t}\|^2 + \frac{\nu^2}{1-\nu}\big(4\beta_{t+1}^2\sum_{i=1}^m\|u^i_{y,t}-\nabla_y f^i(x^i_t,y^i_t)\|^2 \nonumber \\
    & \quad + 4\beta^2_{t+1}m\sigma^2 + 16\eta^2_tL_f^2\sum_{i=1}^m\big( \|\tilde{x}^i_{t+1}-x^i_t\|^2 + \|\tilde{y}^i_{t+1}-y^i_t\|^2\big) \big).\nonumber
\end{align}

\end{lemma}

\begin{proof}
For notational simplicity, let $x_t=[(x^1_t)^T,\cdots,(x^m_t)^T]^T\in \mathbb{R}^{md}$, $\tilde{x}_t=[(\tilde{x}^1_t)^T,\cdots,(\tilde{x}^m_t)^T]^T\in \mathbb{R}^{md}$, $y_t=[(y^1_t)^T,\cdots,(y^m_t)^T]^T\in \mathbb{R}^{mp}$, $\tilde{y}_t=[(\tilde{y}^1_t)^T,\cdots,(\tilde{y}^m_t)^T]^T\in \mathbb{R}^{mp}$, $w_{x,t}=[(w^1_{x,t})^T,\cdots,(w^m_{x,t})^T]^T\in \mathbb{R}^{md}$ and $w_{y,t}=[(w^1_{y,t})^T,\cdots,(w^m_{y,t})^T]^T\in \mathbb{R}^{mp}$ for all $t\geq1$.
By using Assumption \ref{ass:1}, since $W\textbf{1}=\textbf{1}$ and $\tilde{W}=W\otimes I_d$, we have $\tilde{W}(\textbf{1}\otimes \bar{x}_t)=\textbf{1}\otimes\bar{x}_t$.
Meanwhile, we have $\textbf{1}^T(x_t - \textbf{1}\otimes\bar{x})=0$ and $\tilde{W}\textbf{1}=\textbf{1}$.
Thus, we have
\begin{align} \label{eq:D1}
  \|\tilde{W}x_t - \textbf{1}\otimes\bar{x}_t\|^2 = \|\tilde{W}(x_t-\textbf{1}\otimes\bar{x}_t)\|^2 \leq \nu^2\|x_t-\textbf{1}\otimes\bar{x}_t\|^2,
\end{align}
where the above inequality holds by $x_t - \textbf{1}\otimes\bar{x}_t$ is orthogonal to $\textbf{1}$ that is the eigenvector corresponding to the largest eigenvalue of $\tilde{W}$, and $\nu$ denotes the second largest eigenvalue of $\tilde{W}$.

Since $\tilde{x}^i_{t+1} = \sum_{j\in \mathcal{N}_i} W_{i,j}x^j_t - \gamma w^i_{x,t}$ for all $i\in [m]$, we have $\tilde{x}_{t+1} = \tilde{W}x_t - \gamma w_{x,t}$ and $\bar{\tilde{x}}_{t+1} = \bar{x}_t - \gamma \bar{w}_{x,t}$.
Since $x_{t+1}=x_t + \eta_t(\tilde{x}_{t+1}-x_t)$ and $\bar{x}_{t+1}=\bar{x}_t + \eta_t(\bar{\tilde{x}}_{t+1}-\bar{x}_t)$, we have
\begin{align}
  \sum_{i=1}^m\|x^i_{t+1} - \bar{x}_{t+1}\|^2 & = \big\| x_{t+1} - \textbf{1}\otimes\bar{x}_{t+1}\big\|^2 \\
  & = \big\|x_t + \eta_t(\tilde{x}_{t+1}-x_t) - \textbf{1}\otimes\big((\bar{x}_t + \eta_t(\bar{\tilde{x}}_{t+1}-\bar{x}_t)\big) \big\|^2 \nonumber \\
  & \leq (1+\alpha_1)(1-\eta_t)^2\|x_t-\textbf{1}\otimes\bar{x}_t\|^2+(1+\frac{1}{\alpha_1})\eta^2_t
  \|\tilde{x}_{t+1}-\textbf{1}\otimes\bar{\tilde{x}}_{t+1} \|^2 \nonumber \\
  & \mathop{=}^{(i)} (1-\eta_t)\|x_t-\textbf{1}\otimes\bar{x}_t\|^2+\eta_t
  \|\tilde{x}_{t+1}-\textbf{1}\otimes\bar{\tilde{x}}_{t+1} \|^2 \nonumber \\
  & = (1-\eta_t)\|x_t-\textbf{1}\otimes\bar{x}_t\|^2+\eta_t
  \|\tilde{W}x_t - \gamma w_{x,t} -\textbf{1}\otimes\big( \bar{x}_t - \gamma \bar{w}_{x,t}\big) \|^2 \nonumber \\
  & \leq (1-\eta_t)\|x_t-\textbf{1}\otimes\bar{x}_t\|^2+(1+\alpha_2)\eta_t
  \|\tilde{W}x_t - \textbf{1}\otimes\bar{x}_t\|^2 + (1+\frac{1}{\alpha_2})\eta_t\gamma^2\| w_{x,t}-\textbf{1}\otimes\bar{w}_{x,t} \|^2 \nonumber \\
  & \mathop{\leq}^{(ii)} (1-\eta_t)\|x_t-\textbf{1}\otimes\bar{x}_t\|^2+\frac{(1+\nu^2)\eta_t}{2}
  \|x_t - \textbf{1}\otimes\bar{x}_t\|^2 + \frac{\eta_t\gamma^2(1+\nu^2)}{1-\nu^2}\| w_{x,t}-\textbf{1}\otimes\bar{w}_{x,t} \|^2 \nonumber \\
  & \mathop{\leq}^{(iii)} (1-\frac{(1-\nu^2)\eta_t}{2})\|x_t-\textbf{1}\otimes\bar{x}_t\|^2 + \frac{2\eta_t\gamma^2}{1-\nu^2}\| w_{x,t}-\textbf{1}\otimes\bar{w}_{x,t} \|^2 \nonumber \\
  & = (1-\frac{(1-\nu^2)\eta_t}{2})\sum_{i=1}^m\|x^i_t-\bar{x}_t\|^2 + \frac{2\eta_t\gamma^2}{1-\nu^2}\sum_{i=1}^m\| w^i_{x,t}-\bar{w}_{x,t} \|^2,
\end{align}
where the above equality $(i)$ is due to $\alpha_1=\frac{\eta_t}{1-\eta_t}$, and the second inequality $(ii)$ holds by $\alpha_2=\frac{1-\nu^2}{2\nu^2}$ and $\|\tilde{W}x_t - \textbf{1}\otimes\bar{x}_t\|^2 \leq
\nu^2\|x_t - \textbf{1}\otimes\bar{x}_t\|^2$, and
 the above inequality $(ii)$ is due to $0<\nu<1$.
Meanwhile, we have
\begin{align}
 \sum_{i=1}^m\|\tilde{x}^i_{t+1}-\bar{x}_t\|^2 & =\|\tilde{x}_{t+1}-\textbf{1}\otimes\bar{x}_t\|^2
 \nonumber \\
 & = \|\tilde{W}x_t - \gamma w_{x,t} - \textbf{1}\otimes\bar{x}_t\|^2 \nonumber \\
 & \leq (1+\alpha_2)\nu^2\|x_t -\textbf{1}\otimes\bar{x}_t\|^2 + (1+\frac{1}{\alpha_2})\gamma^2\|w_{x,t}\|^2 \nonumber \\
 & \mathop{=}^{(i)} \frac{1+\nu^2}{2}\sum_{i=1}^m\|x^i_t -\bar{x}_t\|^2 + \frac{1+\nu^2}{1-\nu^2}\gamma^2\sum_{i=1}^m\|w^i_{x,t}\|^2,
\end{align}
where the last equality $(i)$ holds by $\alpha_2=\frac{1-\nu^2}{2\nu^2}$.
Then we have
\begin{align}
 \sum_{i=1}^m\|\tilde{x}^i_{t+1}-x^i_t\|^2 & =\|\tilde{x}_{t+1}-x_t\|^2 \nonumber \\
 & = \|\tilde{x}_{t+1}- \textbf{1}\otimes\bar{x}_t + \textbf{1}\otimes\bar{x}_t - x_t\|^2 \nonumber \\
 & \leq 2\|\tilde{x}_{t+1} - \textbf{1}\otimes\bar{x}_t\|^2 + 2\|x_t - \textbf{1}\otimes\bar{x}_t\|^2 \nonumber \\
 & = (3+\nu^2)\sum_{i=1}^m\|x^i_t -\bar{x}_t\|^2 + \frac{2(1+\nu^2)}{1-\nu^2}\gamma^2\sum_{i=1}^m\|w^i_{x,t}\|^2.
\end{align}

Let $w_{x,t}=[(w^1_{x,t})^T,(w^2_{x,t})^T,\cdots,(w^m_{x,t})^T]^T$,
$u_{x,t}=[(u^1_{x,t})^T,(u^2_{x,t})^T,\cdots,(u^m_{x,t})^T]^T$ and $\bar{w}_{x,t} = \frac{1}{m}\sum_{i=1}^mw^i_{x,t}$ and $\bar{u}_{x,t} = \frac{1}{m}\sum_{i=1}^mu^i_{x,t}$. Then we have for any $t\geq 1$,
\begin{align}
  w_{x,t+1} = \tilde{W}\big(w_{x,t} + u_{x,t+1} - u_{x,t}\big). \nonumber
\end{align}
According to the above proof of Lemma~\ref{lem:B1}, we have $\bar{w}_{x,t+1}=\bar{w}_{x,t} + \bar{u}_{x,t+1} - \bar{u}_{x,t}$ for all $t\geq1$.
Thus we have
\begin{align} \label{eq:D3}
     \sum_{i=1}^m\|w^i_{x,t+1}-\bar{w}_{x,t+1}\|^2 & =
     \|w_{x,t+1}-\textbf{1}\otimes\bar{w}_{x,t+1}\|^2 \nonumber \\
    & = \big\| \tilde{W}\big(w_{x,t} + u_{x,t+1} - u_{x,t}\big) -\textbf{1}\otimes\big(\bar{w}_{x,t} + \bar{u}_{x,t+1} - \bar{u}_{x,t} \big)\big\|^2 \nonumber \\
    & \leq (1+c)\|\tilde{W}w_{x,t}-\textbf{1}\otimes\bar{w}_{x,t}\|^2 + (1+\frac{1}{c})\big\|\tilde{W}\big( u_{x,t+1} - u_{x,t}\big)-\textbf{1}\otimes\big(\bar{u}_{x,t+1} - \bar{u}_{x,t}\big)\big\|^2 \nonumber \\
    & \leq (1+c)\nu^2\|w_{x,t}-\textbf{1}\otimes\bar{w}_{x,t}\|^2 + (1+\frac{1}{c})\nu^2\big\|u_{x,t+1} - u_{x,t}-\textbf{1}\otimes\big(\bar{u}_{x,t+1} - \bar{u}_{x,t}\big)\big\|^2 \nonumber \\
    & \leq (1+c)\nu^2\|w_{x,t}-\textbf{1}\otimes\bar{w}_{x,t}\|^2 + (1+\frac{1}{c})\nu^2\big\|u_{x,t+1} - u_{x,t}\big\|^2,
\end{align}
where the first inequality is due to Assumption~\ref{ass:1} as in the above inequality~(\ref{eq:D1}), and the last inequality holds by the inequality $\mathbb{E}\|\xi-\mathbb{E}[\xi]\|^2 \leq \mathbb{E}\|\xi\|^2$.

Since $u^i_{x,t+1} = \nabla f^i(x^i_{t+1};\xi^i_{t+1}) + (1-\alpha_{t+1})(u^i_{x,t}-\nabla f^i(x^i_t;\xi^i_{t+1}))$ for any $i\in [m]$ and $t\geq 1$, we have
\begin{align} \label{eq:D4}
 &\big\|u_{x,t+1} - u_{x,t}\big\|^2 = \sum_{i=1}^m\|u^i_{x,t+1} - u^i_{x,t}\|^2 \nonumber \\
 &= \sum_{i=1}^m\|\nabla_x f^i(x^i_{t+1},y^i_{t+1};\xi^i_{t+1}) - \alpha_{t+1}u^i_{x,t} - (1-\alpha_{t+1})\nabla_x f^i(x^i_t,y^i_t;\xi^i_{t+1})\|^2 \nonumber \\
 &= \sum_{i=1}^m\|\alpha_{t+1}(\nabla_x f^i(x^i_{t+1},y^i_{t+1};\xi^i_{t+1})-\nabla_x f^i(x^i_{t+1},y^i_{t+1})) + \alpha_{t+1}(\nabla_x f^i(x^i_{t+1},y^i_{t+1})-\nabla_x f^i(x^i_t,y^i_t))  \nonumber \\
 & \qquad + \alpha_{t+1}(\nabla_x f^i(x^i_t,y^i_t)-u^i_{x,t}) + (1-\alpha_{t+1})(\nabla_x f^i(x^i_{t+1},y^i_{t+1};\xi^i_{t+1})-\nabla_x f^i(x^i_t,y^i_t;\xi^i_{t+1}))\|^2 \nonumber \\
 &= 4\alpha^2_{t+1}\sum_{i=1}^m\|\nabla_x f^i(x^i_{t+1},y^i_{t+1};\xi^i_{t+1})-\nabla_x f^i(x^i_{t+1},y^i_{t+1})\|^2 +4\alpha^2_{t+1}\sum_{i=1}^m\|\nabla_x f^i(x^i_{t+1},y^i_{t+1})-\nabla_x f^i(x^i_t,y^i_t)\|^2 \nonumber \\
 & \quad + 4\alpha^2_{t+1}\sum_{i=1}^m\|u^i_{x,t}-\nabla_x f^i(x^i_t,y^i_t) \|^2 + 4(1-\alpha_{t+1})^2\sum_{i=1}^m\|\nabla_x f^i(x^i_{t+1},y^i_{t+1};\xi^i_{t+1})-\nabla_x f^i(x^i_t,y^i_t;\xi^i_{t+1})\|^2 \nonumber \\
 & \leq 4\alpha_{t+1}^2\sum_{i=1}^m\|u^i_{x,t}-\nabla_x f^i(x^i_t,y^i_t)\|^2+ 4\alpha^2_{t+1}m\sigma^2 + 8\alpha_{t+1}^2L_f^2\sum_{i=1}^m\big( \|\tilde{x}^i_{t+1}-x^i_t\|^2 + \|\tilde{y}^i_{t+1}-y^i_t\|^2\big) \nonumber \\
 & \quad + 8(1-\alpha_{t+1})^2L_f^2\sum_{i=1}^m\big( \|\tilde{x}^i_{t+1}-x^i_t\|^2 + \|\tilde{y}^i_{t+1}-y^i_t\|^2\big) \nonumber \\
 & \leq 4\alpha_{t+1}^2\sum_{i=1}^m\|u^i_{x,t}-\nabla_x f^i(x^i_t,y^i_t)\|^2+ 4\alpha^2_{t+1}m\sigma^2 + 16\eta^2_tL_f^2\sum_{i=1}^m\big( \|\tilde{x}^i_{t+1}-x^i_t\|^2 + \|\tilde{y}^i_{t+1}-y^i_t\|^2\big),
\end{align}
where the last inequality holds by $0<\alpha_t<1$ and $x^i_{t+1}=x^i_t+\eta_t(\tilde{x}^i_{t+1}-x^i_t)$.

Plugging the above inequalities (\ref{eq:D4}) into (\ref{eq:D3}), we have
\begin{align}
    \sum_{i=1}^m\|w^i_{x,t+1}-\bar{w}_{x,t+1}\|^2
    & \leq (1+c)\nu^2\|w_{x,t}-\textbf{1}\otimes\bar{w}_{x,t}\|^2 + (1+\frac{1}{c})\nu^2\big\|u_{x,t+1} - u_{x,t}\big\|^2 \nonumber \\
    & \leq (1+c)\nu^2\|w_{x,t}-\textbf{1}\otimes\bar{w}_{x,t}\|^2 + (1+\frac{1}{c})\nu^2\big(4\alpha_{t+1}^2\sum_{i=1}^m\|u^i_{x,t}-\nabla_x f^i(x^i_t,y^i_t)\|^2+ 4\alpha^2_{t+1}m\sigma^2 \nonumber \\
    & \quad + 16\eta^2_tL_f^2\sum_{i=1}^m\big( \|\tilde{x}^i_{t+1}-x^i_t\|^2 + \|\tilde{y}^i_{t+1}-y^i_t\|^2\big) \big).
\end{align}
Let $c=\frac{1}{\nu}-1$, we have
\begin{align}
   \sum_{i=1}^m\|w^i_{x,t+1}-\bar{w}_{x,t+1}\|^2
    & \leq \nu\sum_{i=1}^m\|w^i_{x,t}-\bar{w}_{x,t}\|^2 + \frac{\nu^2}{1-\nu}\big(4\alpha_{t+1}^2\sum_{i=1}^m\|u^i_{x,t}-\nabla_x f^i(x^i_t,y^i_t)\|^2 \nonumber \\
    & \quad + 4\alpha^2_{t+1}m\sigma^2 + 16\eta^2_tL_f^2\sum_{i=1}^m\big( \|\tilde{x}^i_{t+1}-x^i_t\|^2 + \|\tilde{y}^i_{t+1}-y^i_t\|^2\big) \big).
\end{align}

Similarly, we can get
\begin{align}
\sum_{i=1}^m\|y^i_{t+1} - \bar{y}_{t+1}\|^2 & \leq (1-\frac{(1-\nu^2)\eta_t}{2})\sum_{i=1}^m\|y^i_t-\bar{y}_t\|^2
+ \frac{2\eta_t\lambda^2}{1-\nu^2}\sum_{i=1}^m\| w^i_{y,t}-\bar{w}_{y,t} \|^2,
\nonumber \\
\sum_{i=1}^m\|\tilde{y}^i_{t+1}-y^i_t\|^2
 & \leq (3+\nu^2)\sum_{i=1}^m\|y^i_t -\bar{y}_t\|^2 + \frac{2(1+\nu^2)}{1-\nu^2}\lambda^2
 \sum_{i=1}^m\|w^i_{y,t}\|^2, \nonumber \\
\sum_{i=1}^m\|w^i_{y,t+1}-\bar{w}_{y,t+1}\|^2
    & \leq \nu\sum_{i=1}^m\|w^i_{y,t}-\bar{w}_{y,t}\|^2 + \frac{\nu^2}{1-\nu}\big(4\beta_{t+1}^2\sum_{i=1}^m\|u^i_{y,t}-\nabla_y f^i(x^i_t,y^i_t)\|^2 \nonumber \\
    & \quad + 4\beta^2_{t+1}m\sigma^2 + 16\eta^2_tL_f^2\sum_{i=1}^m\big( \|\tilde{x}^i_{t+1}-x^i_t\|^2 + \|\tilde{y}^i_{t+1}-y^i_t\|^2\big) \big).\nonumber
\end{align}

\end{proof}

\begin{lemma} \label{lem:E1}
 Suppose the sequence $\{\bar{x}_t\}_{t=1}^T$ be generated from Algorithm~\ref{alg:1}. Given $ 0< \gamma \leq \frac{1}{4L\eta_t}$ for all $t\geq 0$,
 then we have
 \begin{align}
 F(\bar{x}_{t+1}) \leq F(\bar{x}_t)+\frac{2\gamma\eta_t}{\rho}\|\nabla F(\bar{x}_t)-\bar{u}_{x,t}\|^2+\frac{\gamma\eta_t}{2}\frac{1}{m}\sum_{i=1}^m\|w^i_{x,t} -\bar{w}_{x,t}\|^2-\frac{\gamma\eta_t}{4} \frac{1}{m}\sum_{i=1}^m\|w^i_{x,t}\|^2.
 \end{align}
\end{lemma}

\begin{proof}
Since $\tilde{x}^i_{t+1} = \sum_{j\in \mathcal{N}_i} W_{i,j}x^j_t - \gamma w^i_{x,t}$, we have
\begin{align}
\bar{\tilde{x}}_{t+1} & = \frac{1}{m}\sum_{i=1}^m\tilde{x}^i_{t+1} =\frac{1}{m}\sum_{i=1}^m\sum_{j\in \mathcal{N}_i} W_{i,j}x^j_t - \gamma \frac{1}{m}\sum_{i=1}^mw^i_{x,t} \nonumber \\
& = \frac{1}{m}\sum_{j\in \mathcal{N}_i}x^j_t\sum_{i=1}^m W_{i,j} - \gamma \frac{1}{m}\sum_{i=1}^mw^i_{x,t} =
\bar{x}_t - \gamma \bar{w}_{x,t},
\end{align}
where the last equality is due to $\sum_{i=1}^m W_{i,j}=1$.
Since $\bar{u}_{x,t}=\bar{w}_{x,t}$, we have
\begin{align}
 \langle w^i_{x,t}, w^i_{x,t}\rangle \mathop{=}^{(i)} \langle w^i_{x,t} -\bar{w}_{x,t} , w^i_{x,t}\rangle + \langle\bar{u}_{x,t} , w^i_{x,t}\rangle  \leq \frac{1}{2}\|w^i_{x,t} -\bar{w}_{x,t}\|^2 + \frac{1}{2}\|w^i_{x,t}\|^2 + \langle\bar{u}_{x,t} , w^i_{x,t}\rangle,
\end{align}
we have
\begin{align}
 0 \leq  \frac{1}{2}\|w^i_{x,t} -\bar{w}_{x,t}\|^2 - \frac{1}{2}\|w^i_{x,t}\|^2 + \langle\bar{u}_{x,t}, w^i_{x,t}\rangle.
\end{align}
Then taking an average over $i$ from 1 to $m$ yields that
\begin{align} \label{eq:E1}
 0 & \leq  \frac{\gamma\eta_t}{2}\frac{1}{m}\sum_{i=1}^m\|w^i_{x,t} -\bar{w}_{x,t}\|^2 - \frac{\gamma\eta_t}{2}\frac{1}{m}\sum_{i=1}^m\|w^i_{x,t}\|^2 + \gamma\eta_t\frac{1}{m}\sum_{i=1}^m\langle\bar{u}_{x,t}, w^i_{x,t}\rangle \nonumber \\
 & = \frac{\gamma\eta_t}{2}\frac{1}{m}\sum_{i=1}^m\|w^i_{x,t} -\bar{w}_{x,t}\|^2 - \frac{\gamma\eta_t}{2}\frac{1}{m}\sum_{i=1}^m\|w^i_{x,t}\|^2 + \gamma\eta_t\langle\bar{u}_{x,t}, \bar{w}_{x,t}\rangle.
\end{align}

According to Assumption~\ref{ass:1}, i.e., the function $F(x)$ is $L$-smooth,
we have
 \begin{align} \label{eq:E2}
  F(\bar{x}_{t+1}) & \leq F(\bar{x}_t) + \langle\nabla F(\bar{x}_t), \bar{x}_{t+1}-\bar{x}_t\rangle + \frac{L}{2}\|\bar{x}_{t+1}-\bar{x}_t\|^2 \\
  & = F(\bar{x}_t)+ \eta_t\langle \nabla F(\bar{x}_t),\bar{\tilde{x}}_{t+1}-x_t\rangle + \frac{L\eta_t^2}{2}\|\bar{\tilde{x}}_{t+1}-\bar{x}_t\|^2 \nonumber \\
   & = F(\bar{x}_t)+ \eta_t\langle \nabla F(\bar{x}_t)-\bar{u}_{x,t} + \bar{u}_{x,t},-\gamma\bar{w}_{x,t}\rangle + \frac{L\gamma^2\eta_t^2}{2}\|\bar{w}_{x,t}\|^2 \nonumber \\
  & \leq F(\bar{x}_t) + 2\gamma\eta_t\|\nabla F(\bar{x}_t)-\bar{u}_{x,t}\|^2+\frac{\gamma\eta_t}{8}\|\bar{w}_{x,t}\|^2 -\gamma \eta_t\langle \bar{u}_{x,t},\bar{w}_{x,t}\rangle + \frac{L\gamma^2\eta_t^2}{2}\|\bar{w}_{x,t}\|^2 \nonumber \\
  & \leq  F(\bar{x}_t)+ 2\gamma\eta_t\|\nabla F(\bar{x}_t)-\bar{u}_{x,t}\|^2+\frac{\gamma\eta_t}{8}\frac{1}{m}\sum_{i=1}^m\|w^i_{x,t}\|^2 -\gamma \eta_t\langle \bar{u}_{x,t},\bar{w}_{x,t}\rangle + \frac{L\gamma^2\eta_t^2}{2}\frac{1}{m}\sum_{i=1}^m\|w^i_{x,t}\|^2, \nonumber
 \end{align}
where the second equality is due to $\bar{x}_{t+1}=\bar{x}_t + \eta_t(\bar{\tilde{x}}_{t+1}-\bar{x}_t)$.

By summing the above inequalities (\ref{eq:E1}) with (\ref{eq:E2}), we can obtain
 \begin{align}
   F(\bar{x}_{t+1}) & \leq  F(\bar{x}_t) + 2\gamma\eta_t\|\nabla F(\bar{x}_t)-\bar{u}_{x,t}\|^2+\frac{\gamma\eta_t}{2}\frac{1}{m}\sum_{i=1}^m\|w^i_{x,t} -\bar{w}_{x,t}\|^2-(\frac{3\gamma\eta_t}{8} - \frac{L\gamma^2\eta_t^2}{2})\frac{1}{m}\sum_{i=1}^m\|w^i_{x,t}\|^2 \nonumber \\
   & \leq F(\bar{x}_t)+\frac{2\gamma\eta_t}{\rho}\|\nabla F(\bar{x}_t)-\bar{u}_{x,t}\|^2+\frac{\gamma\eta_t}{2}\frac{1}{m}\sum_{i=1}^m\|w^i_{x,t} -\bar{w}_{x,t}\|^2-\frac{\gamma\eta_t}{4} \frac{1}{m}\sum_{i=1}^m\|w^i_{x,t}\|^2,
 \end{align}
where the last inequality is due to $\gamma \leq \frac{1}{4L\eta_t}$ for all $t\geq 1$.

\end{proof}

\begin{lemma} \label{lem:F1}
(Restatement of Lemma 3)
Suppose the sequence $\{\bar{x}_t,\bar{y}_t\}_{t=1}^T$ be generated from Algorithm~\ref{alg:1}.
Under the Assumptions~\ref{ass:2}-\ref{ass:3}, given $0<\gamma\leq \frac{\lambda\mu}{16L}$ and $0<\lambda \leq \frac{1}{2L_f\eta_t}$ for all $t\geq 1$, we have
\begin{align}
F(\bar{x}_{t+1}) - f(\bar{x}_{t+1},\bar{y}_{t+1})
& \leq (1-\frac{\eta_t\lambda\mu}{2}) \big(F(\bar{x}_t) -f(\bar{x}_t,\bar{y}_t)\big) + \frac{\eta_t\gamma}{8}\|\bar{w}_{x,t}\|^2  -\frac{\eta_t\lambda}{4}\|\bar{w}_{y,t}\|^2 \nonumber \\
& \quad + \eta_t\lambda\|\nabla_y f(\bar{x}_t,\bar{y}_t)-\bar{w}_{y,t}\|^2,
\end{align}
where $F(\bar{x}_t)=f(\bar{x}_t,y^*(\bar{x}_t))$ with $y^*(\bar{x}_t) \in \arg\max_{y}f(\bar{x}_t,y)$ for all $t\geq 1$.
\end{lemma}

\begin{proof}
Using $L_f$-smoothness of $f(x,\cdot)$, such that
\begin{align}
    f(\bar{x}_{t+1},\bar{y}_t) + \langle \nabla_y f(\bar{x}_{t+1},\bar{y}_t), \bar{y}_{t+1}-\bar{y}_t \rangle - \frac{L_f}{2}\|\bar{y}_{t+1}-\bar{y}_t\|^2 \leq f(\bar{x}_{t+1},\bar{y}_{t+1}),
\end{align}
then we have
\begin{align} \label{eq:F1}
    f(\bar{x}_{t+1},\bar{y}_t) & \leq f(\bar{x}_{t+1},\bar{y}_{t+1}) - \langle \nabla_y f(\bar{x}_{t+1},\bar{y}_t), \bar{y}_{t+1}-\bar{y}_t \rangle + \frac{L_f}{2}\|\bar{y}_{t+1}-\bar{y}_t\|^2 \nonumber \\
    & = f(\bar{x}_{t+1},\bar{y}_{t+1}) - \eta_t\langle \nabla_y f(\bar{x}_{t+1},\bar{y}_t), \bar{\tilde{y}}_{t+1}-\bar{y}_t \rangle + \frac{L_f\eta^2_t}{2}\|\bar{\tilde{y}}_{t+1}-\bar{y}_t\|^2.
\end{align}

Next, we bound the inner product in \eqref{eq:F1}. According to the line 6 of Algorithm~\ref{alg:1}, i.e.,
$\bar{\tilde{y}}_{t+1} = \bar{y}_t + \lambda \bar{w}_{y,t}$, we have
\begin{align} \label{eq:F2}
    & - \eta_t\langle \nabla_y f(\bar{x}_{t+1},\bar{y}_t), \bar{\tilde{y}}_{t+1}-\bar{y}_t \rangle \nonumber \\
    & = - \eta_t\lambda\langle \nabla_y f(\bar{x}_{t+1},\bar{y}_t), \bar{w}_{y,t} \rangle \nonumber \\
    & = -\frac{\eta_t\lambda}{2}\Big( \|\nabla_y f(\bar{x}_{t+1},\bar{y}_t)\|^2 + \|\bar{w}_{y,t}\|^2 - \|\nabla_y f(\bar{x}_{t+1},\bar{y}_t)-\nabla_y f(\bar{x}_t,\bar{y}_t) + \nabla_y f(\bar{x}_t,\bar{y}_t)-\bar{w}_{y,t}\|^2 \Big) \nonumber \\
    & \leq -\frac{\eta_t\lambda}{2} \|\nabla_y f(\bar{x}_{t+1},\bar{y}_t)\|^2 -\frac{\eta_t}{2\lambda} \|\bar{\tilde{y}}_{t+1}-\bar{y}_t\|^2 + \eta_t\lambda L^2_f \|\bar{x}_{t+1}-\bar{x}_t\|^2 + \eta_t\lambda\|\nabla_y f(\bar{x}_t,\bar{y}_t)-\bar{w}_{y,t}\|^2 \nonumber \\
    & \leq -\eta_t\lambda\mu\big(F(\bar{x}_{t+1})-f(\bar{x}_{t+1},\bar{y}_t)\big)-\frac{\eta_t}{2\lambda} \|\bar{\tilde{y}}_{t+1}-\bar{y}_t\|^2 + \eta_t\lambda L^2_f \|\bar{x}_{t+1}-\bar{x}_t\|^2 + \eta_t\lambda\|\nabla_y f(\bar{x}_t,\bar{y}_t)-\bar{w}_{y,t}\|^2,
\end{align}
where the last inequality is due to the quadratic growth condition of $\mu$-PL functions, i.e.,
\begin{align}
    \|\nabla_y f(\bar{x}_{t+1},\bar{y}_t)\|^2 \geq 2\mu\big( \max_{y'}f(\bar{x}_{t+1},y')-f(\bar{x}_{t+1},\bar{y}_t)\big) = 2\mu\big( F(\bar{x}_{t+1})-f(\bar{x}_{t+1},\bar{y}_t)\big).
\end{align}
Substituting \eqref{eq:F2} in \eqref{eq:F1}, we have
\begin{align} \label{eq:F3}
    f(\bar{x}_{t+1},\bar{y}_t)
    & = f(\bar{x}_{t+1},\bar{y}_{t+1})-\eta_t\lambda\mu\big(F(\bar{x}_{t+1})-f(\bar{x}_{t+1},\bar{y}_t)\big)-\frac{\eta_t}{2\lambda} \|\bar{\tilde{y}}_{t+1}-\bar{y}_t\|^2 + \eta_t\lambda L^2_f \|\bar{x}_{t+1}-\bar{x}_t\|^2 \nonumber \\
    & \quad + \eta_t\lambda\|\nabla_y f(\bar{x}_t,\bar{y}_t)-\bar{w}_{y,t}\|^2 + \frac{L_f\eta^2_t}{2}\|\bar{\tilde{y}}_{t+1}-\bar{y}_t\|^2,
\end{align}
then rearranging the terms, we can obtain
\begin{align} \label{eq:F4}
    F(\bar{x}_{t+1}) - f(\bar{x}_{t+1},\bar{y}_{t+1})
    & = (1-\eta_t\lambda\mu)\big(F(\bar{x}_{t+1})-f(\bar{x}_{t+1},\bar{y}_t)\big)-\frac{\eta_t}{2\lambda} \|\bar{\tilde{y}}_{t+1}-\bar{y}_t\|^2 + \eta_t\lambda L^2_f \|\bar{x}_{t+1}-\bar{x}_t\|^2 \nonumber \\
    & \quad + \eta_t\lambda\|\nabla_y f(\bar{x}_t,\bar{y}_t)-\bar{w}_{y,t}\|^2 + \frac{L_f\eta^2_t}{2}\|\bar{\tilde{y}}_{t+1}-\bar{y}_t\|^2.
\end{align}

Next, using $L_f$-smoothness of function $f(\cdot,\bar{y}_t)$, such that
\begin{align}
    f(\bar{x}_t,\bar{y}_t) + \langle \nabla_x f(\bar{x}_t,\bar{y}_t), \bar{x}_{t+1}-\bar{x}_t \rangle - \frac{L_f}{2}\|\bar{x}_{t+1}-\bar{x}_t\|^2 \leq f(\bar{x}_{t+1},\bar{y}_t),
\end{align}
then we have
\begin{align}
    & f(\bar{x}_t,\bar{y}_t) -f(\bar{x}_{t+1},\bar{y}_t)  \nonumber \\
    & \leq -\langle \nabla_x f(\bar{x}_t,\bar{y}_t), \bar{x}_{t+1}-\bar{x}_t \rangle + \frac{L_f}{2}\|\bar{x}_{t+1}-\bar{x}_t\|^2 \nonumber \\
    & = -\eta_t\langle \nabla_x f(\bar{x}_t,\bar{y}_t) - \nabla F(\bar{x}_t), \bar{\tilde{x}}_{t+1}-\bar{x}_t \rangle - \eta_t\langle \nabla F(\bar{x}_t), \bar{\tilde{x}}_{t+1}-\bar{x}_t \rangle + \frac{L_f\eta^2_t}{2}\|\bar{\tilde{x}}_{t+1}-\bar{x}_t\|^2 \nonumber \\
    & \leq \frac{\eta_t}{8\gamma}\|\bar{\tilde{x}}_{t+1}-\bar{x}_t\|^2 + 2\eta_t\gamma\|\nabla_x f(\bar{x}_t,\bar{y}_t) - \nabla F(\bar{x}_t)\|^2  - \eta_t\langle \nabla F(\bar{x}_t), \bar{\tilde{x}}_{t+1}-\bar{x}_t \rangle + \frac{L_f\eta^2_t}{2}\|\bar{\tilde{x}}_{t+1}-\bar{x}_t\|^2 \nonumber \\
    & \leq \frac{\eta_t}{8\gamma}\|\bar{\tilde{x}}_{t+1}-\bar{x}_t\|^2 + 2L^2_f\eta_t\gamma \|\bar{y}_t - y^*(\bar{x}_t)\|^2 + F(\bar{x}_t) - F(\bar{x}_{t+1}) \nonumber \\
    & \quad  + \frac{\eta^2_tL}{2}\|\bar{\tilde{x}}_{t+1}-\bar{x}_t\|^2 + \frac{\eta^2_tL_f}{2}\|\bar{\tilde{x}}_{t+1}-\bar{x}_t\|^2 \nonumber \\
    & \leq \frac{4L^2_f\eta_t\gamma}{\mu} \big(F(\bar{x}_t) -f(\bar{x}_t,\bar{y}_t)\big) + F(\bar{x}_t) - F(\bar{x}_{t+1}) + \eta_t(\frac{1}{8\gamma}+\eta_tL)\|\bar{\tilde{x}}_{t+1}-\bar{x}_t\|^2,
\end{align}
where the second last inequality is due to Lemma~\ref{lem:A1}, i.e.,
$L$-smoothness of function $F(x)$, and the last inequality holds by Lemma~\ref{lem:A2} and $L_f\leq L$.
Then we have
\begin{align} \label{eq:F5}
    F(\bar{x}_{t+1}) - f(\bar{x}_{t+1},\bar{y}_t) & = F(\bar{x}_{t+1}) - F(\bar{x}_t) + F(\bar{x}_t)- f(\bar{x}_t,\bar{y}_t) + f(\bar{x}_t,\bar{y}_t) -f(\bar{x}_{t+1},\bar{y}_t) \nonumber \\
    & \leq (1+\frac{4L^2_f\eta_t\gamma}{\mu}) \big(F(\bar{x}_t) -f(\bar{x}_t,\bar{y}_t)\big) + \eta_t(\frac{1}{8\gamma}+\eta_tL)\|\bar{\tilde{x}}_{t+1}-\bar{x}_t\|^2.
\end{align}

Substituting \eqref{eq:F5} in \eqref{eq:F4}, we get
\begin{align}
    & F(\bar{x}_{t+1}) - f(\bar{x}_{t+1},\bar{y}_{t+1}) \nonumber \\
    & \leq (1-\eta_t\lambda\mu)(1+\frac{4L^2_f\eta_t\gamma}{\mu}) \big(F(\bar{x}_t) -f(\bar{x}_t,\bar{y}_t)\big) + \eta_t(\frac{1}{8\gamma}+\eta_tL)(1-\eta_t\lambda\mu)\|\bar{\tilde{x}}_{t+1}-\bar{x}_t\|^2 \nonumber \\
    & \quad -\frac{\eta_t}{2\lambda} \|\bar{\tilde{y}}_{t+1}-\bar{y}_t\|^2 + \eta_t\lambda L^2_f \|\bar{x}_{t+1}-\bar{x}_t\|^2 + \eta_t\lambda\|\nabla_y f(\bar{x}_t,\bar{y}_t)-\bar{w}_{y,t}\|^2 + \frac{L_f\eta^2_t}{2}\|\bar{\tilde{y}}_{t+1}-\bar{y}_t\|^2 \nonumber \\
    & = (1-\eta_t\lambda\mu)(1+\frac{4L^2_f\eta_t\gamma}{\mu}) \big(F(\bar{x}_t) -f(\bar{x}_t,\bar{y}_t)\big) + \eta_t\big(\frac{1}{8\gamma}+\eta_tL-\frac{\eta_t\lambda\mu}{8\gamma}-\eta^2_tL\lambda\mu+\eta^2_tL^2_f\lambda\big)\|\bar{\tilde{x}}_{t+1}-\bar{x}_t\|^2 \nonumber \\
    & \quad -\frac{\eta_t}{2}\big(\frac{1}{\lambda}-L_f\eta_t\big) \|\bar{\tilde{y}}_{t+1}-\bar{y}_t\|^2 + \eta_t\lambda\|\nabla_y f(\bar{x}_t,\bar{y}_t)-\bar{w}_{y,t}\|^2 \nonumber \\
    & \leq (1-\frac{\eta_t\lambda\mu}{2}) \big(F(\bar{x}_t) -f(\bar{x}_t,\bar{y}_t)\big) + \frac{\eta_t}{8\gamma}\|\bar{\tilde{x}}_{t+1}-\bar{x}_t\|^2  -\frac{\eta_t}{4\lambda}\|\bar{\tilde{y}}_{t+1}-\bar{y}_t\|^2 + \eta_t\lambda\|\nabla_y f(\bar{x}_t,\bar{y}_t)-\bar{w}_{y,t}\|^2,
\end{align}
where the last inequality holds by $L=L_f(1+\frac{\kappa}{2})$,  $\gamma\leq \frac{\lambda\mu}{16L}$ and $\lambda\leq \frac{1}{2L_f\eta_t}$ for all $t\geq 1$, i.e.,
\begin{align}
   & \gamma\leq \frac{\lambda\mu}{16L} \Rightarrow \lambda \geq \frac{16L\gamma}{\mu} =  16(\kappa+\frac{\kappa^2}{2})\gamma\geq 8\kappa^2\gamma \Rightarrow  \frac{\eta_t\lambda\mu}{2} \geq \frac{4L^2_f\eta_t\gamma}{\mu} \nonumber \\
   & L=L_f(1+\kappa) \Rightarrow \eta^2_tL\lambda\mu \geq \eta^2_tL^2_f\lambda \nonumber \\
   &\lambda \leq \frac{1}{2\eta_tL_f}  \Rightarrow \frac{1}{2\lambda} \geq \eta_t L_f, \ \forall t\geq 1.
\end{align}

Since $\bar{\tilde{x}}_{t+1}= \bar{x}_t - \gamma \bar{w}_{x,t}$ and $\bar{\tilde{y}}_{t+1}=\bar{y}_t + \lambda \bar{w}_{y,t}$, we can get
\begin{align}
    & F(\bar{x}_{t+1}) - f(\bar{x}_{t+1},\bar{y}_{t+1}) \nonumber \\
    & \leq (1-\frac{\eta_t\lambda\mu}{2}) \big(F(\bar{x}_t) -f(\bar{x}_t,\bar{y}_t)\big) + \frac{\eta_t}{8\gamma}\|\bar{\tilde{x}}_{t+1}-\bar{x}_t\|^2  -\frac{\eta_t}{4\lambda}\|\bar{\tilde{y}}_{t+1}-\bar{y}_t\|^2 + \eta_t\lambda\|\nabla_y f(\bar{x}_t,\bar{y}_t)-\bar{w}_{y,t}\|^2  \nonumber \\
    & \leq (1-\frac{\eta_t\lambda\mu}{2}) \big(F(\bar{x}_t) -f(\bar{x}_t,\bar{y}_t)\big) + \frac{\eta_t\gamma}{8}\|\bar{w}_{x,t}\|^2  -\frac{\eta_t\lambda}{4}\|\bar{w}_{y,t}\|^2 + \eta_t\lambda\|\nabla_y f(\bar{x}_t,\bar{y}_t)-\bar{w}_{y,t}\|^2.
\end{align}

\end{proof}

\begin{theorem}  \label{th:A1}
(Restatement of Theorem 1)
 Suppose the sequences $\{x_t,y_t\}_{t=1}^T$ be generated from Algorithm~\ref{alg:1}.
 Under the above Assumptions~\ref{ass:1}-\ref{ass:4}, and let $\alpha_t =\beta_t = O(\frac{1}{T^{2/3}})$ and $\eta_t=\eta=O(\frac{1}{T^{1/3}})$ for all $t\geq1$, $0<\gamma\leq \frac{\lambda\mu}{16L}$ and $\lambda =O(1)$ for all $t\geq 0$, we have
\begin{align}
  \frac{1}{T}\sum_{t=1}^T \|\nabla F(\bar{x}_t)\|
  \leq O(\frac{1}{T^{1/3}}+\frac{\sigma^2}{T^{1/3}}).
\end{align}
\end{theorem}

\begin{proof}
Without loss of generality, let $\eta=\eta_1=\cdots=\eta_T$.
According to Lemma \ref{lem:E1}, we have
\begin{align} \label{eq:H1}
 F(\bar{x}_{t+1}) \leq F(\bar{x}_t)+2\gamma\eta\|\nabla F(\bar{x}_t)-\bar{u}_{x,t}\|^2+\frac{\gamma\eta}{2}\frac{1}{m}\sum_{i=1}^m\|w^i_{x,t} -\bar{w}_{x,t}\|^2-\frac{\gamma\eta}{4} \frac{1}{m}\sum_{i=1}^m\|w^i_{x,t}\|^2.
\end{align}

According to the Lemma~\ref{lem:C1},
we have
 \begin{align} \label{eq:H2}
 \mathbb{E}\|\bar{u}_{x,t} - \overline{\nabla_x f(x_{t},y_{t})}\|^2
  & \leq (1-\alpha_{t})\mathbb{E} \|\bar{u}_{x,t-1} -\overline{\nabla_x f(x_{t-1},y_{t-1})}\|^2 + \frac{2\alpha_{t}^2\sigma^2}{m} \nonumber \\
  & \quad + \frac{2L^2_f\eta^2}{m^2}\sum_{i=1}^m\mathbb{E}\big(\|\tilde{x}^i_{t}-x^i_{t-1}\|^2+\|\tilde{y}^i_{t}-y^i_{t-1}\|^2\big),
  \end{align}
and
 \begin{align} \label{eq:H3}
 \frac{1}{m}\sum_{i=1}^m\mathbb{E}\|u^i_{x,t} - \nabla_x f^i(x^i_t,y^i_t)\|^2
  & \leq (1-\alpha_{t})\frac{1}{m}\sum_{i=1}^m\mathbb{E} \|u^i_{x,t-1} - \nabla_x f^i(x^i_{t-1},y^i_{t-1})\|^2 + 2\alpha_{t}^2\sigma^2 \nonumber \\
  & \quad + 2L^2_f\eta^2\frac{1}{m}\sum_{i=1}^m\mathbb{E}\big(\|\tilde{x}^i_{t}-x^i_{t-1}\|^2+\|\tilde{y}^i_{t}-y^i_{t-1}\|^2\big).
 \end{align}

Similarly, we also have
 \begin{align} \label{eq:H4}
 \mathbb{E}\|\bar{u}_{y,t} - \overline{\nabla_y f(x_{t},y_{t})}\|^2
 & \leq (1-\beta_{t})\mathbb{E} \|\bar{u}_{y,t-1} -\overline{\nabla_y f(x_{t-1},y_{t-1})}\|^2 + \frac{2\beta_{t}^2\sigma^2}{m} \nonumber \\
 & \quad + \frac{2L^2_f\eta^2}{m^2}\sum_{i=1}^m\mathbb{E}\big(\|\tilde{x}^i_{t}-x^i_{t-1}\|^2+\|\tilde{y}^i_{t}-y^i_{t-1}\|^2\big),
 \end{align}
 and
 \begin{align} \label{eq:H5}
 \frac{1}{m}\sum_{i=1}^m\mathbb{E}\|u^i_{y,t} - \nabla_y f^i(x^i_t,y^i_t)\|^2
  & \leq (1-\beta_{t})\frac{1}{m}\sum_{i=1}^m\mathbb{E} \|u^i_{y,t-1} - \nabla_y f^i(x^i_{t-1},y^i_{t-1})\|^2 + 2\beta_{t}^2\sigma^2 \nonumber \\
  & \quad + 2L^2_f\eta^2\frac{1}{m}\sum_{i=1}^m\mathbb{E}\big(\|\tilde{x}^i_{t}-x^i_{t-1}\|^2+\|\tilde{y}^i_{t}-y^i_{t-1}\|^2\big).
 \end{align}

According to Lemma \ref{lem:D1}, we have
\begin{align} \label{eq:H6}
\frac{1}{m}\sum_{i=1}^m\|w^i_{x,t}-\bar{w}_{x,t}\|^2
    & \leq \nu\frac{1}{m}\sum_{i=1}^m\|w^i_{x,t-1}-\bar{w}_{x,t-1}\|^2 + \frac{\nu^2}{1-\nu}\big(4\alpha_{t}^2\frac{1}{m}\sum_{i=1}^m\|u^i_{x,t-1}-\nabla_x f^i(x^i_{t-1},y^i_{t-1})\|^2
    \nonumber \\
    & \quad + 4\alpha^2_{t}\sigma^2 + 16\eta^2L^2_f\frac{1}{m}\sum_{i=1}^m\big( \|\tilde{x}^i_{t}-x^i_{t-1}\|^2 + \|\tilde{y}^i_{t}-y^i_{t-1}\|^2\big) \big),
\end{align}
and
\begin{align} \label{eq:H7}
\frac{1}{m}\sum_{i=1}^m\|x^i_t - \bar{x}_t\|^2 & \leq (1-\frac{(1-\nu^2)\eta}{2})\frac{1}{m}\sum_{i=1}^m\|x^i_{t-1}-\bar{x}_{t-1}\|^2
+ \frac{2\eta\gamma^2}{1-\nu^2}\frac{1}{m}\sum_{i=1}^m\| w^i_{x,t-1}-\bar{w}_{x,t-1} \|^2 \nonumber \\
& \leq (1-\frac{(1-\nu^2)\eta}{2})\frac{1}{m}\sum_{i=1}^m\|x^i_{t-1}-\bar{x}_{t-1}\|^2
+ \frac{2\eta\gamma^2}{1-\nu^2}\frac{1}{m}\sum_{i=1}^m(\| w^i_{x,t-1}\|^2 + \|\bar{w}_{x,t-1} \|^2) \nonumber \\
& \leq (1-\frac{(1-\nu^2)\eta}{2})\frac{1}{m}\sum_{i=1}^m\|x^i_{t-1}-\bar{x}_{t-1}\|^2
+ \frac{4\eta\gamma^2}{1-\nu^2}\frac{1}{m}\sum_{i=1}^m\|w^i_{x,t-1}\|^2,
\end{align}
and
\begin{align} \label{eq:H8}
\frac{1}{m}\sum_{i=1}^m\|\tilde{x}^i_t-x^i_{t-1}\|^2
 & \leq (3+\nu^2)\frac{1}{m}\sum_{i=1}^m\|x^i_{t-1} -\bar{x}_{t-1}\|^2 + \frac{2(1+\nu^2)}{1-\nu^2}\gamma^2
 \frac{1}{m}\sum_{i=1}^m\|w^i_{x,t-1}\|^2.
\end{align}
Similarly, we can get
\begin{align} \label{eq:H9}
\frac{1}{m}\sum_{i=1}^m\|w^i_{y,t}-\bar{w}_{y,t}\|^2
    & \leq \nu\frac{1}{m}\sum_{i=1}^m\|w^i_{y,t-1}-\bar{w}_{y,t-1}\|^2 + \frac{\nu^2}{1-\nu}\big(4\beta_{t}^2\frac{1}{m}\sum_{i=1}^m\|u^i_{y,t-1}-\nabla_y f^i(x^i_{t-1},y^i_{t-1})\|^2
    \nonumber \\
    & \quad + 4\beta^2_{t}\sigma^2 + 16\eta^2L^2_f\frac{1}{m}\sum_{i=1}^m\big( \|\tilde{x}^i_{t}-x^i_{t-1}\|^2 + \|\tilde{y}^i_{t}-y^i_{t-1}\|^2\big) \big),
\end{align}
and
\begin{align} \label{eq:H10}
\frac{1}{m}\sum_{i=1}^m\|y^i_t - \bar{y}_t\|^2 & \leq (1-\frac{(1-\nu^2)\eta}{2})\frac{1}{m}\sum_{i=1}^m\|y^i_{t-1}-\bar{y}_{t-1}\|^2
+ \frac{4\eta\lambda^2}{1-\nu^2}\frac{1}{m}\sum_{i=1}^m\| w^i_{y,t-1}\|^2,
\end{align}
and
\begin{align} \label{eq:H11}
\frac{1}{m}\sum_{i=1}^m\|\tilde{y}^i_{t}-y^i_{t-1}\|^2
 & \leq (3+\nu^2)\frac{1}{m}\sum_{i=1}^m\|y^i_{t-1} -\bar{y}_{t-1}\|^2 + \frac{2(1+\nu^2)}{1-\nu^2}\lambda^2
 \frac{1}{m}\sum_{i=1}^m\|w^i_{y,t-1}\|^2.
\end{align}

Since $\nabla F(\bar{x}_t)=\frac{1}{m}\sum_{i=1}^m\nabla_x f^i(\bar{x}_t,y^*(\bar{x}_t))$, we have
\begin{align}
\|\bar{u}_{x,t}-\nabla F(\bar{x}_t)\|^2 & = \|\bar{u}_{x,t}- \overline{\nabla_x f(x_t,y_t)} + \overline{\nabla_x f(x_t,y_t)} - \nabla F(\bar{x}_t)\|^2  \\
& \leq 2\|\bar{u}_{x,t}- \overline{\nabla_x f(x_t,y_t)}\|^2 + 2\|\overline{\nabla_x f(x_t,y_t)} - \nabla F(\bar{x}_t)\|^2 \nonumber \\
& \leq 2\|\bar{u}_{x,t}- \overline{\nabla_x f(x_t,y_t)}\|^2 + 2\|\frac{1}{m}\sum_{i=1}^m\nabla_x f^i(x^i_t,y^i_t) - \frac{1}{m}\sum_{i=1}^m\nabla_x f^i(\bar{x}_t,y^*(\bar{x}_t))\|^2 \nonumber \\
& \leq 2\|\bar{u}_{x,t}- \overline{\nabla_x f(x_t,y_t)}\|^2 + \frac{4L^2_f}{m}\sum_{i=1}^m\big(\|x^i_t - \bar{x}_t\|^2+\|y^i_t - y^*(\bar{x}_t)\|^2\big) \nonumber \\
& \leq 2\|\bar{u}_{x,t}- \overline{\nabla_x f(x_t,y_t)}\|^2 + \frac{4L^2_f}{m}\sum_{i=1}^m\big(\|x^i_t - \bar{x}_t\|^2 + 2\|y^i_t - \bar{y}_t\|^2 + 2\|\bar{y}_t - y^*(\bar{x}_t)\|^2\big) \nonumber \\
& \leq 2\|\bar{u}_{x,t}- \overline{\nabla_x f(x_t,y_t)}\|^2 + \frac{4L^2_f}{m}\sum_{i=1}^m\big(\|x^i_t - \bar{x}_t\|^2 + 2\|y^i_t - \bar{y}_t\|^2\big) + \frac{16L^2_f}{\mu}\big(F(\bar{x}_t) - f(\bar{x}_t,\bar{y}_t)\big), \nonumber
\end{align}
where the last inequality is due to Assumption~\ref{ass:2}. Then we can obtain
\begin{align} \label{eq:H12}
 - \|\bar{u}_{x,t}- \overline{\nabla_x f(x_t,y_t)}\|^2 \leq -\frac{1}{2}\|\bar{u}_{x,t} - \nabla F(\bar{x}_t)\|^2 + \frac{2L^2_f}{m}\sum_{i=1}^m\big(\|x^i_t - \bar{x}_t\|^2+2\|y^i_t - \bar{y}_t\|^2\big) +
  \frac{8L^2_f}{\mu}\big(F(\bar{x}_t) - f(\bar{x}_t,\bar{y}_t)\big).
\end{align}
Since $\bar{u}_{x,t}=\bar{w}_{x,t}$ for all $t\geq1$, we have
\begin{align}
 \frac{1}{m}\sum_{i=1}^m\|w^i_{x,t}-\nabla F(x^i_t)\|^2 & = \frac{1}{m}\sum_{i=1}^m\|w^i_{x,t}-\bar{w}_{x,t}+\bar{u}_{x,t}-\nabla F(\bar{x}_t)+\nabla F(\bar{x}_t)-\nabla F(x^i_t)\|^2 \nonumber \\
 & \leq 3\frac{1}{m}\sum_{i=1}^m\|w^i_{x,t}-\bar{w}_{x,t}\|^2 + 3\|\bar{u}_{x,t}-\nabla F(\bar{x}_t)\|^2+3\frac{1}{m}\sum_{i=1}^m\|\nabla F(\bar{x}_t)-\nabla F(x^i_t)\|^2 \nonumber \\
 & \leq 3\frac{1}{m}\sum_{i=1}^m\|w^i_{x,t}-\bar{w}_{x,t}\|^2 + 3\|\bar{u}_{x,t}-\nabla F(\bar{x}_t)\|^2+3L^2\frac{1}{m}\sum_{i=1}^m\|x^i_t-\bar{x}_t\|^2. 
\end{align}
Then we have
\begin{align} \label{eq:H13}
 -\|\bar{u}_{x,t}-\nabla F(\bar{x}_t)\|^2
 & \leq -\frac{1}{3m}\sum_{i=1}^m\|w^i_{x,t}-\nabla F(x^i_t)\|^2 + \frac{1}{m}\sum_{i=1}^m\|w^i_{x,t}-\bar{w}_{x,t}\|^2 +\frac{L^2}{m}\sum_{i=1}^m\|x^i_t-\bar{x}_t\|^2.
\end{align}
According to Lemma~\ref{lem:F1}, we have
\begin{align} \label{eq:H14}
F(\bar{x}_{t+1}) - f(\bar{x}_{t+1},\bar{y}_{t+1})
& \leq (1-\frac{\eta\lambda\mu}{2}) \big(F(\bar{x}_t) -f(\bar{x}_t,\bar{y}_t)\big) + \frac{\eta\gamma}{8}\|\bar{w}_{x,t}\|^2  -\frac{\eta\lambda}{4}\|\bar{w}_{y,t}\|^2 + \eta\lambda\|\nabla_y f(\bar{x}_t,\bar{y}_t)-\bar{w}_{y,t}\|^2 \nonumber \\
& \leq (1-\frac{\eta\lambda\mu}{2}) \big(F(\bar{x}_t) -f(\bar{x}_t,\bar{y}_t)\big) + \frac{\eta\gamma}{8}\|\bar{w}_{x,t}\|^2  -\frac{\eta\lambda}{4}\|\bar{w}_{y,t}\|^2 \nonumber \\
& \quad + 2\eta\lambda\|\nabla_y f(\bar{x}_t,\bar{y}_t)-\overline{\nabla_y f(x_t,y_t)}\|^2 +2\eta\lambda\|\overline{\nabla_y f(x_t,y_t)}-\bar{w}_{y,t}\|^2 \nonumber \\
& \leq (1-\frac{\eta\lambda\mu}{2}) \big(F(\bar{x}_t) -f(\bar{x}_t,\bar{y}_t)\big) + \frac{\eta\gamma}{8}\|\bar{w}_{x,t}\|^2  -\frac{\eta\lambda}{4}\|\bar{w}_{y,t}\|^2 \nonumber \\
& \quad + \frac{2\eta\lambda}{m}\sum_{i=1}^m\|\nabla_y f^i(\bar{x}_t,\bar{y}_t)-\nabla_y f^i(x^i_t,y^i_t)\|^2 +2\eta\lambda\|\overline{\nabla_y f(x_t,y_t)}-\bar{w}_{y,t}\|^2 \nonumber \\
& \leq (1-\frac{\eta\lambda\mu}{2}) \big(F(\bar{x}_t) -f(\bar{x}_t,\bar{y}_t)\big) + \frac{\eta\gamma}{8m}\sum_{i=1}^m\|w^i_{x,t}\|^2  -\frac{\eta\lambda}{4}\|\bar{w}_{y,t}\|^2 \nonumber \\
& \quad + \frac{4\eta\lambda L^2_f}{m}\sum_{i=1}^m\big( \|\bar{x}_t-x^i_t\|^2 + \|\bar{y}_t-y^i_t\|^2 \big) + 2\eta\lambda\|\overline{\nabla_y f(x_t,y_t)}-\bar{w}_{y,t}\|^2,
\end{align}
where the last inequality is due to Assumption~\ref{ass:2}.

Next, we define a useful Lyapunov function (i.e., potential function), for any $t\geq 1$
\begin{align} \label{eq:H15}
\Omega_t & = \mathbb{E}_t \Big[ F(\bar{x}_t) + \frac{72\gamma L^2_f}{\lambda\mu^2}\big(F(\bar{x}_t) -f(\bar{x}_t,\bar{y}_t)\big) + (\rho_{x,t-1}-\frac{9\gamma\eta}{2})\|\bar{u}_{x,t-1} - \overline{\nabla_x f(x_{t-1},y_{t-1})}\|^2  \nonumber \\
& \qquad + (\rho_{y,t-1}-\frac{144\gamma\eta L^2_f}{\mu^2})\|\bar{u}_{y,t-1} - \overline{\nabla_y f(x_{t-1},y_{t-1})}\|^2 + \varrho_{x,t-1}\frac{1}{m}\sum_{i=1}^m \|u^i_{x,t-1} - \nabla_x f^i(x^i_{t-1},y^i_{t-1})\|^2   \nonumber \\
& \qquad + \varrho_{y,t-1}\frac{1}{m}\sum_{i=1}^m \|u^i_{y,t-1} - \nabla_y f^i(x^i_{t-1},y^i_{t-1})\|^2 + (\theta_{x,t-1}-9\gamma\eta L_f^2-\frac{\gamma\eta L^2}{2}-\frac{288\gamma\eta L^4_f}{\mu^2})\frac{1}{m}\sum_{i=1}^m\|x^i_{t-1}-\bar{x}_{t-1}\|^2  \nonumber \\
& \qquad  + (\theta_{y,t-1}-18\gamma\eta L_f^2-\frac{288\gamma\eta L^4_f}{\mu^2})\frac{1}{m}\sum_{i=1}^m\|y^i_{t-1}-\bar{y}_{t-1}\|^2 +(\vartheta_{x,t-1}-\frac{3\gamma\eta}{4})\frac{1}{m}\sum_{i=1}^m\|w^i_{x,t-1}-\bar{w}_{x,t-1}\|^2
 \nonumber \\
& \qquad +\vartheta_{y,t-1}\frac{1}{m}\sum_{i=1}^m\|w^i_{y,t-1}-\bar{w}_{y,t-1}\|^2 +
\frac{\gamma\eta}{12}\frac{1}{m}\sum_{i=1}^m\|w^i_{x,t-1}\|^2 + \frac{18\gamma L^2_f\eta}{\mu^2}\|\bar{w}_{y,t-1}\|^2\Big],
\end{align}
where $\rho_{x,t-1}\geq \frac{9\gamma\eta}{2}$, $\rho_{y,t-1}\geq \frac{144\gamma\eta L^2_f}{\mu^2}$, $\varrho_{x,t-1}\geq 0$, $\varrho_{y,t-1}\geq 0$,
$\theta_{x,t-1}\geq 9\gamma\eta L_f^2+\frac{\gamma\eta L^2}{2}+\frac{288\gamma\eta L^4_f}{\mu^2}$, $\theta_{y,t-1}\geq 18\gamma\eta L_f^2+\frac{288\gamma\eta L^4_f}{\mu^2}$,
$\vartheta_{x,t-1}\geq \frac{3\gamma\eta}{4}$ and $\vartheta_{y,t-1}\geq 0$ for all $t\geq1$.

Then we have
\begin{align} \label{eq:H16}
 \Omega_{t+1} & = \mathbb{E}_{t+1} \Big[ F(\bar{x}_{t+1}) + \frac{72\gamma L^2_f}{\lambda\mu^2}\big(F(\bar{x}_{t+1}) - f(\bar{x}_{t+1},\bar{y}_{t+1})\big) + (\rho_{x,t}-\frac{9\gamma\eta}{2})\|\bar{u}_{x,t} - \overline{\nabla_x f(x_{t},y_{t})}\|^2 \nonumber \\
& \qquad + (\rho_{y,t}-\frac{144\gamma\eta L^2_f}{\mu^2})\|\bar{u}_{y,t} - \overline{\nabla_y f(x_{t},y_{t})}\|^2 + \varrho_{x,t}\frac{1}{m}\sum_{i=1}^m \|u^i_{x,t} - \nabla_x f^i(x^i_{t},y^i_{t})\|^2   + \varrho_{y,t}\frac{1}{m}\sum_{i=1}^m \|u^i_{y,t} - \nabla_y f^i(x^i_{t},y^i_{t})\|^2  \nonumber \\
& \qquad + (\theta_{x,t}-9\gamma\eta L_f^2-\frac{\gamma\eta L^2}{2}-\frac{288\gamma\eta L^4_f}{\mu^2})\frac{1}{m}\sum_{i=1}^m\|x^i_{t}-\bar{x}_{t}\|^2  + (\theta_{y,t}-18\gamma\eta L_f^2-\frac{288\gamma\eta L^4_f}{\mu^2})\frac{1}{m}\sum_{i=1}^m\|y^i_{t}-\bar{y}_{t}\|^2 \nonumber \\
& \qquad +(\vartheta_{x,t}-\frac{3\gamma\eta}{4})\frac{1}{m}\sum_{i=1}^m\|w^i_{x,t}-\bar{w}_{x,t}\|^2
+\vartheta_{y,t}\frac{1}{m}\sum_{i=1}^m\|w^i_{y,t}-\bar{w}_{y,t}\|^2 +
\frac{\gamma\eta}{12}\frac{1}{m}\sum_{i=1}^m\|w^i_{x,t}\|^2 + \frac{18\gamma L^2_f\eta}{\mu^2}\|\bar{w}_{y,t}\|^2\Big] \nonumber \\
 & \mathop{\leq}^{(i)} \mathbb{E}_{t+1}\Big[ F(\bar{x}_t) + \frac{72\gamma L^2_f}{\lambda\mu^2}\big(F(\bar{x}_t) -f(\bar{x}_t,\bar{y}_t)\big) -\frac{\gamma\eta}{4}\|\bar{u}_{x,t} - \nabla F(\bar{x}_t)\|^2  - (\frac{\gamma\eta}{6} - \frac{9\gamma^2\eta L^2_f}{\lambda\mu^2})\frac{1}{m}\sum_{i=1}^m\|w^i_{x,t}\|^2 \nonumber \\
 & \qquad  + \frac{9\gamma\eta L^2_f}{m}
 \sum_{i=1}^m\big( \|x^i_t-\bar{x}_t\|^2 + 2\|y^i_t-\bar{y}_t\|^2 \big)+\frac{\gamma\eta}{2}\frac{1}{m}\sum_{i=1}^m\|w^i_{x,t} -\bar{w}_{x,t}\|^2  \nonumber \\
 & \qquad  + \frac{288\gamma\eta L^4_f}{\mu^2m}\sum_{i=1}^m\big( \|\bar{x}_t-x^i_t\|^2 + \|\bar{y}_t-y^i_t\|^2 \big) + \rho_{x,t}\|\bar{u}_{x,t}-\overline{\nabla_x f(x_t,y_t)}\|^2 + \rho_{y,t}\|\bar{u}_{y,t} - \overline{\nabla_y f(x_{t},y_{t})}\|^2\nonumber \\
& \qquad + \varrho_{x,t}\frac{1}{m}\sum_{i=1}^m \|u^i_{x,t} - \nabla_x f^i(x^i_{t},y^i_{t})\|^2   + \varrho_{y,t}\frac{1}{m}\sum_{i=1}^m \|u^i_{y,t} - \nabla_y f^i(x^i_{t},y^i_{t})\|^2  \nonumber \\
& \qquad + (\theta_{x,t}-9\gamma\eta L_f^2-\frac{\gamma\eta L^2}{2}-\frac{288\gamma\eta L^4_f}{\mu^2})\frac{1}{m}\sum_{i=1}^m\|x^i_{t}-\bar{x}_{t}\|^2  + (\theta_{y,t}-18\gamma\eta L_f^2-\frac{288\gamma\eta L^4_f}{\mu^2})\frac{1}{m}\sum_{i=1}^m\|y^i_{t}-\bar{y}_{t}\|^2 \nonumber \\
& \qquad +(\vartheta_{x,t}-\frac{3\gamma\eta}{4})\frac{1}{m}\sum_{i=1}^m\|w^i_{x,t}-\bar{w}_{x,t}\|^2
+\vartheta_{y,t}\frac{1}{m}\sum_{i=1}^m\|w^i_{y,t}-\bar{w}_{y,t}\|^2 \Big] \nonumber \\
 & \mathop{\leq}^{(ii)} \mathbb{E}_{t+1}\Big[ F(\bar{x}_t) + \frac{72\gamma L^2_f}{\lambda\mu^2}\big(F(\bar{x}_t) -f(\bar{x}_t,\bar{y}_t)\big) -\frac{\gamma\eta}{12m}\sum_{i=1}^m\|w^i_{x,t}-\nabla F(x^i_t)\|^2  - (\frac{\gamma\eta}{6} - \frac{9\gamma^2\eta L^2_f}{\lambda\mu^2})\frac{1}{m}\sum_{i=1}^m\|w^i_{x,t}\|^2 \nonumber \\
 & \qquad  + \frac{9\gamma\eta L^2_f}{m}
 \sum_{i=1}^m\big( \|x^i_t-\bar{x}_t\|^2 + 2\|y^i_t-\bar{y}_t\|^2 \big)+\frac{3\gamma\eta}{4}\frac{1}{m}\sum_{i=1}^m\|w^i_{x,t} -\bar{w}_{x,t}\|^2 +\frac{\gamma\eta L^2}{4m}\sum_{i=1}^m\|x^i_t-\bar{x}_t\|^2 \nonumber \\
 & \qquad + \frac{288\gamma\eta L^4_f}{\mu^2m}\sum_{i=1}^m\big( \|\bar{x}_t-x^i_t\|^2 + \|\bar{y}_t-y^i_t\|^2 \big) + \rho_{x,t}\|\bar{u}_{x,t}-\overline{\nabla_x f(x_t,y_t)}\|^2 \nonumber \\
& \qquad + \rho_{y,t}\|\bar{u}_{y,t} - \overline{\nabla_y f(x_{t},y_{t})}\|^2+ \varrho_{x,t}\frac{1}{m}\sum_{i=1}^m \|u^i_{x,t} - \nabla_x f^i(x^i_{t},y^i_{t})\|^2   + \varrho_{y,t}\frac{1}{m}\sum_{i=1}^m \|u^i_{y,t} - \nabla_y f^i(x^i_{t},y^i_{t})\|^2  \nonumber \\
& \qquad + (\theta_{x,t}-9\gamma\eta L_f^2-\frac{\gamma\eta L^2}{2}-\frac{288\gamma\eta L^4_f}{\mu^2})\frac{1}{m}\sum_{i=1}^m\|x^i_{t}-\bar{x}_{t}\|^2  + (\theta_{y,t}-18\gamma\eta L_f^2-\frac{288\gamma\eta L^4_f}{\mu^2})\frac{1}{m}\sum_{i=1}^m\|y^i_{t}-\bar{y}_{t}\|^2 \nonumber \\
& \qquad +(\vartheta_{x,t}-\frac{3\gamma\eta}{4})\frac{1}{m}\sum_{i=1}^m\|w^i_{x,t}-\bar{w}_{x,t}\|^2
+\vartheta_{y,t}\frac{1}{m}\sum_{i=1}^m\|w^i_{y,t}-\bar{w}_{y,t}\|^2 \Big],
\end{align}
where the inequality (i) is due to the above inequalities (\ref{eq:H1}) and (\ref{eq:H14}); and
the inequality (ii) holds by the above inequalities (\ref{eq:H13}).

Then we have
\begin{align} \label{eq:H17}
\Omega_{t+1} & \leq  \mathbb{E}_{t+1}\Big[ F(\bar{x}_t) + \frac{72\gamma L^2_f}{\lambda\mu^2}\big(F(\bar{x}_t) -f(\bar{x}_t,\bar{y}_t)\big) -\frac{\gamma\eta}{12m}\sum_{i=1}^m\|w^i_{x,t}-\nabla F(x^i_t)\|^2 -\frac{\gamma\eta L^2}{4}\frac{1}{m}\sum_{i=1}^m\|x^i_{t}-\bar{x}_{t}\|^2 \nonumber \\
 & \qquad - (\frac{\gamma\eta}{6} - \frac{9\gamma^2\eta L^2_f}{\lambda\mu^2})\frac{1}{m}\sum_{i=1}^m\|w^i_{x,t}\|^2 + \rho_{x,t}\|\bar{u}_{x,t}-\overline{\nabla_x f(x_t,y_t)}\|^2 + \rho_{y,t}\|\bar{u}_{y,t} - \overline{\nabla_y f(x_{t},y_{t})}\|^2
 \nonumber \\
& \qquad + \varrho_{x,t}\frac{1}{m}\sum_{i=1}^m \|u^i_{x,t} - \nabla_x f^i(x^i_{t},y^i_{t})\|^2   + \varrho_{y,t}\frac{1}{m}\sum_{i=1}^m \|u^i_{y,t} - \nabla_y f^i(x^i_{t},y^i_{t})\|^2 + \theta_{x,t}\frac{1}{m}\sum_{i=1}^m\|x^i_{t}-\bar{x}_{t}\|^2  \nonumber \\
& \qquad  + \theta_{y,t}\frac{1}{m}\sum_{i=1}^m\|y^i_{t}-\bar{y}_{t}\|^2  +\vartheta_{x,t}\frac{1}{m}\sum_{i=1}^m\|w^i_{x,t}-\bar{w}_{x,t}\|^2
+\vartheta_{y,t}\frac{1}{m}\sum_{i=1}^m\|w^i_{y,t}-\bar{w}_{y,t}\|^2 \Big] \nonumber \\
& \mathop{\leq}^{(i)} \mathbb{E}_{t+1}\Big[ F(\bar{x}_t) + \frac{72\gamma L^2_f}{\lambda\mu^2}\big(F(\bar{x}_t) -f(\bar{x}_t,\bar{y}_t)\big) -\frac{\gamma\eta}{12m}\sum_{i=1}^m\|w^i_{x,t}-\nabla F(x^i_t)\|^2-\frac{\gamma\eta L^2}{4}\frac{1}{m}\sum_{i=1}^m\|x^i_{t}-\bar{x}_{t}\|^2 \nonumber \\
 & \qquad  - (\frac{\gamma\eta}{6} - \frac{9\gamma^2\eta L^2_f}{\lambda\mu^2})\frac{1}{m}\sum_{i=1}^m\|w^i_{x,t}\|^2  \nonumber \\
& \qquad + \rho_{x,t}(1-\alpha_{t})\mathbb{E} \|\bar{u}_{x,t-1} -\overline{\nabla_x f(x_{t-1},y_{t-1})}\|^2 + \frac{2\rho_{x,t}\alpha_{t}^2\sigma^2}{m} + \frac{2\rho_{x,t}L^2_f\eta^2}{m^2}\sum_{i=1}^m\mathbb{E}\big(\|\tilde{x}^i_{t}-x^i_{t-1}\|^2+\|\tilde{y}^i_{t}-y^i_{t-1}\|^2\big) \nonumber \\
& \qquad + \rho_{y,t}(1-\beta_{t})\mathbb{E} \|\bar{u}_{y,t-1} -\overline{\nabla_y f(x_{t-1},y_{t-1})}\|^2 + \frac{2\rho_{y,t}\beta_{t}^2\sigma^2}{m} + \frac{2\rho_{y,t}L^2_f\eta^2}{m^2}\sum_{i=1}^m\mathbb{E}\big(\|\tilde{x}^i_{t}-x^i_{t-1}\|^2+\|\tilde{y}^i_{t}-y^i_{t-1}\|^2\big) \nonumber \\
& \qquad + \varrho_{x,t}(1-\alpha_{t})\frac{1}{m}\sum_{i=1}^m\mathbb{E} \|u^i_{x,t-1} - \nabla_x f^i(x^i_{t-1},y^i_{t-1})\|^2 + 2\varrho_{x,t}\alpha_{t}^2\sigma^2 + \frac{2\varrho_{x,t}L^2_f\eta^2}{m}\sum_{i=1}^m\mathbb{E}\big(\|\tilde{x}^i_{t}-x^i_{t-1}\|^2+\|\tilde{y}^i_{t}-y^i_{t-1}\|^2\big)
  \nonumber \\
&\qquad + \varrho_{y,t}(1-\beta_{t})\frac{1}{m}\sum_{i=1}^m\mathbb{E} \|u^i_{y,t-1} - \nabla_y f^i(x^i_{t-1},y^i_{t-1})\|^2 + 2\varrho_{y,t}\beta_{t}^2\sigma^2 + \frac{2\varrho_{y,t}L^2_f\eta^2}{m}\sum_{i=1}^m\mathbb{E}\big(\|\tilde{x}^i_{t}-x^i_{t-1}\|^2+\|\tilde{y}^i_{t}-y^i_{t-1}\|^2\big)  \nonumber \\
& \qquad + \theta_{x,t}(1-\frac{(1-\nu^2)\eta}{2})\frac{1}{m}\sum_{i=1}^m\|x^i_{t-1}-\bar{x}_{t-1}\|^2
+ \frac{4\theta_{x,t}\eta\gamma^2}{1-\nu^2}\frac{1}{m}\sum_{i=1}^m\|w^i_{x,t-1}\|^2 \nonumber \\
& \qquad + \theta_{y,t}(1-\frac{(1-\nu^2)\eta}{2})\frac{1}{m}\sum_{i=1}^m\|y^i_{t-1}-\bar{y}_{t-1}\|^2
+ \frac{4\theta_{y,t}\eta\lambda^2}{1-\nu^2}\frac{1}{m}\sum_{i=1}^m\| w^i_{y,t-1}\|^2 \nonumber \\
& \qquad +\vartheta_{x,t}\nu\frac{1}{m}\sum_{i=1}^m\|w^i_{x,t-1}-\bar{w}_{x,t-1}\|^2 + \frac{\vartheta_{x,t}\nu^2}{1-\nu}\Big(4\alpha_{t}^2\frac{1}{m}\sum_{i=1}^m\|u^i_{x,t-1}-\nabla_x f^i(x^i_{t-1},y^i_{t-1})\|^2
+ 4\alpha^2_{t}\sigma^2 \nonumber \\
& \qquad + 16\eta^2L^2_f\frac{1}{m}\sum_{i=1}^m\big( \|\tilde{x}^i_{t}-x^i_{t-1}\|^2 + \|\tilde{y}^i_{t}-y^i_{t-1}\|^2\big) \Big)    \nonumber \\
& \qquad
+\vartheta_{y,t}\nu\frac{1}{m}\sum_{i=1}^m\|w^i_{y,t-1}-\bar{w}_{y,t-1}\|^2 + \frac{\vartheta_{y,t}\nu^2}{1-\nu}\Big( 4\beta_{t}^2\frac{1}{m}\sum_{i=1}^m\|u^i_{y,t-1}-\nabla_y f^i(x^i_{t-1},y^i_{t-1})\|^2 + 4\beta^2_{t}\sigma^2 \nonumber \\
& \qquad + 16\eta^2L^2_f\frac{1}{m}\sum_{i=1}^m\big( \|\tilde{x}^i_{t}-x^i_{t-1}\|^2 + \|\tilde{y}^i_{t}-y^i_{t-1}\|^2\big) \Big) \Big],
\end{align}
where the inequality (i) is due to the above inequalities~(\ref{eq:H2})-(\ref{eq:H6}) and~(\ref{eq:H9}).

Since $0<\alpha_t < 1$ for all $t\geq1$, let $\rho_{x,t} = \frac{9\gamma\eta}{2\alpha_t}\geq \frac{9\gamma\eta}{2}$ and $\rho_{x,t} \leq \rho_{x,t-1}$, then we have $\rho_{x,t}(1-\alpha_t)\leq \rho_{x,t-1} - \frac{9\gamma\eta}{2}$. Since $0<\beta_t < 1$ for all $t\geq1$, let $\rho_{y,t} = \frac{144\gamma\eta L^2_f}{\mu^2\beta_t} \geq \frac{144\gamma\eta L^2_f}{\mu^2}$ and $\rho_{y,t} \leq \rho_{y,t-1}$, then we have $\rho_{y,t}(1-\beta_t)\leq \rho_{y,t-1} - \frac{144\gamma\eta L^2_f}{\mu^2}$.
Since $0<\alpha_t < 1$ for all $t\geq1$, let $\varrho_{x,t} =\frac{4\vartheta_{x,t}\nu^2}{1-\nu}\geq \frac{4\vartheta_{x,t}\alpha_t\nu^2}{1-\nu}$ and $\varrho_{x,t}\leq \varrho_{x,t-1}$, we have $\varrho_{x,t} -\varrho_{x,t}\alpha_t + \frac{4\vartheta_{x,t}\alpha^2_t\nu^2}{1-\nu} \leq \varrho_{x,t-1}$.
Since $0<\beta_t < 1$ for all $t\geq1$, let $\varrho_{y,t} =\frac{4\vartheta_{y,t}\nu^2}{1-\nu}\geq \frac{4\vartheta_{y,t}\beta_t\nu^2}{1-\nu}$ and $\varrho_{y,t}\leq \varrho_{y,t-1}$, we have $\varrho_{y,t} -\varrho_{y,t}\beta_t + \frac{4\vartheta_{y,t}\beta^2_t\nu^2}{1-\nu} \leq \varrho_{y,t-1}$.

Let $\vartheta_{x,t}=\frac{\gamma\eta}{1-\nu}$ for all $t\geq1$, since $0<\nu<1$, we have $\nu\vartheta_{x,t}=\vartheta_{x,t}-(1-\nu)\vartheta_{x,t} \leq \vartheta_{x,t}-\frac{3\gamma\eta}{4}=\vartheta_{x,t-1}-\frac{3\gamma\eta}{4}$.
Meanwhile, let $\Delta_t = \frac{2\rho_{x,t}\alpha_{t}^2\sigma^2}{m}+ \frac{2\rho_{y,t}\beta_{t}^2\sigma^2}{m}+ 2\varrho_{x,t}\alpha_{t}^2\sigma^2 + 2\varrho_{y,t}\beta_{t}^2\sigma^2 + \frac{4\nu^2\vartheta_{x,t} \alpha^2_{t}\sigma^2}{1-\nu} + \frac{4 \nu^2\vartheta_{y,t}\beta^2_{t}\sigma^2}{1-\nu}$ and
$H_t=\frac{2\rho_{x,t}L^2_f\eta^2}{m}+ \frac{2\rho_{y,t}L^2_f\eta^2}{m}+ 2\varrho_{x,t}L^2_f\eta^2+ 2\varrho_{y,t}L^2_f\eta^2+\frac{16\eta^2L^2_f\vartheta_{x,t}\nu^2}{1-\nu}+\frac{16\eta^2L^2_f\vartheta_{y,t}\nu^2}{1-\nu}$.

Based on the choice of these parameters and the above inequality (\ref{eq:H17}), we have
\begin{align} \label{eq:H18}
\Omega_{t+1} & \leq  \mathbb{E}_{t+1}\Big[ F(\bar{x}_t) + \frac{72\gamma L^2_f}{\lambda\mu^2}\big(F(\bar{x}_t) -f(\bar{x}_t,\bar{y}_t)\big) -\frac{\gamma\eta}{12m}\sum_{i=1}^m\|w^i_{x,t}-\nabla F(x^i_t)\|^2 -\frac{\gamma\eta L^2}{4}\frac{1}{m}\sum_{i=1}^m\|x^i_{t}-\bar{x}_{t}\|^2 \nonumber \\
 & \qquad - (\frac{\gamma\eta}{6} - \frac{9\gamma^2\eta L^2_f}{\lambda\mu^2})\frac{1}{m}\sum_{i=1}^m\|w^i_{x,t}\|^2 + \Delta_t + \frac{H_t}{m}\sum_{i=1}^m\mathbb{E}\big(\|\tilde{x}^i_{t}-x^i_{t-1}\|^2+\|\tilde{y}^i_{t}-y^i_{t-1}\|^2\big)  \nonumber \\
& \qquad + (\rho_{x,t-1} - \frac{9\gamma\eta}{2})\mathbb{E} \|\bar{u}_{x,t-1} -\overline{\nabla_x f(x_{t-1},y_{t-1})}\|^2 + (\rho_{y,t-1}- \frac{144\gamma\eta L^2_f}{\mu^2})\mathbb{E} \|\bar{u}_{y,t-1} -\overline{\nabla_y f(x_{t-1},y_{t-1})}\|^2 \nonumber \\
& \qquad  + \varrho_{x,t-1}\frac{1}{m}\sum_{i=1}^m\mathbb{E} \|u^i_{x,t-1} - \nabla_x f^i(x^i_{t-1},y^i_{t-1})\|^2 + \varrho_{y,t-1}\frac{1}{m}\sum_{i=1}^m\mathbb{E} \|u^i_{y,t-1} - \nabla_y f^i(x^i_{t-1},y^i_{t-1})\|^2 \nonumber \\
& \qquad + \theta_{x,t}(1-\frac{(1-\nu^2)\eta}{2})\frac{1}{m}\sum_{i=1}^m\|x^i_{t-1}-\bar{x}_{t-1}\|^2
+ \frac{4\theta_{x,t}\eta\gamma^2}{1-\nu^2}\frac{1}{m}\sum_{i=1}^m\|w^i_{x,t-1}\|^2 \nonumber \\
& \qquad + \theta_{y,t}(1-\frac{(1-\nu^2)\eta}{2})\frac{1}{m}\sum_{i=1}^m\|y^i_{t-1}-\bar{y}_{t-1}\|^2
+ \frac{4\theta_{y,t}\eta\lambda^2}{1-\nu^2}\frac{1}{m}\sum_{i=1}^m\| w^i_{y,t-1}\|^2 \nonumber \\
& \qquad +(\vartheta_{x,t-1}-\frac{3\gamma\eta}{4})\frac{1}{m}\sum_{i=1}^m\|w^i_{x,t-1}-\bar{w}_{x,t-1}\|^2
+\nu\vartheta_{y,t}\frac{1}{m}\sum_{i=1}^m\|w^i_{y,t-1}-\bar{w}_{y,t-1}\|^2 \Big] \nonumber \\
& \mathop{\leq}^{(i)}  \mathbb{E}_{t+1}\Big[ F(\bar{x}_t) + \frac{72\gamma L^2_f}{\lambda\mu^2}\big(F(\bar{x}_t) -f(\bar{x}_t,\bar{y}_t)\big) -\frac{\gamma\eta}{12m}\sum_{i=1}^m\|w^i_{x,t}-\nabla F(x^i_t)\|^2-\frac{\gamma\eta L^2}{4}\frac{1}{m}\sum_{i=1}^m\|x^i_{t}-\bar{x}_{t}\|^2   \nonumber \\
 & \qquad  - (\frac{\gamma\eta}{6} - \frac{9\gamma^2\eta L^2_f}{\lambda\mu^2})\frac{1}{m}\sum_{i=1}^m\|w^i_{x,t}\|^2 + \Delta_t + (3+\nu^2)H_t\frac{1}{m}\sum_{i=1}^m\|x^i_{t-1} -\bar{x}_{t-1}\|^2 + H_t\frac{2(1+\nu^2)}{1-\nu^2}\gamma^2
 \frac{1}{m}\sum_{i=1}^m\|w^i_{x,t-1}\|^2 \nonumber \\
 & \qquad +(3+\nu^2)H_t\frac{1}{m}\sum_{i=1}^m\|y^i_{t-1} -\bar{y}_{t-1}\|^2 + H_t\frac{2(1+\nu^2)}{1-\nu^2}\lambda^2
 \frac{1}{m}\sum_{i=1}^m\|w^i_{y,t-1}\|^2\nonumber \\
& \qquad + (\rho_{x,t-1} - \frac{9\gamma\eta}{2})\mathbb{E} \|\bar{u}_{x,t-1} -\overline{\nabla_x f(x_{t-1},y_{t-1})}\|^2 + (\rho_{y,t-1}- \frac{144\gamma\eta L^2_f}{\mu^2})\mathbb{E} \|\bar{u}_{y,t-1} -\overline{\nabla_y f(x_{t-1},y_{t-1})}\|^2 \nonumber \\
& \qquad  + \varrho_{x,t-1}\frac{1}{m}\sum_{i=1}^m\mathbb{E} \|u^i_{x,t-1} - \nabla_x f^i(x^i_{t-1},y^i_{t-1})\|^2 + \varrho_{y,t-1}\frac{1}{m}\sum_{i=1}^m\mathbb{E} \|u^i_{y,t-1} - \nabla_y f^i(x^i_{t-1},y^i_{t-1})\|^2 \nonumber \\
& \qquad + \theta_{x,t}(1-\frac{(1-\nu^2)\eta}{2})\frac{1}{m}\sum_{i=1}^m\|x^i_{t-1}-\bar{x}_{t-1}\|^2
+ \frac{4\theta_{x,t}\eta\gamma^2}{1-\nu^2}\frac{1}{m}\sum_{i=1}^m\|w^i_{x,t-1}\|^2 \nonumber \\
& \qquad + \theta_{y,t}(1-\frac{(1-\nu^2)\eta}{2})\frac{1}{m}\sum_{i=1}^m\|y^i_{t-1}-\bar{y}_{t-1}\|^2
+ \frac{4\theta_{y,t}\eta\lambda^2}{1-\nu^2}\frac{1}{m}\sum_{i=1}^m\| w^i_{y,t-1}\|^2 \nonumber \\
& \qquad +(\vartheta_{x,t-1}-\frac{3\gamma\eta}{4})\frac{1}{m}\sum_{i=1}^m\|w^i_{x,t-1}-\bar{w}_{x,t-1}\|^2
+\nu\vartheta_{y,t}\frac{1}{m}\sum_{i=1}^m\|w^i_{y,t-1}-\bar{w}_{y,t-1}\|^2 \Big],
\end{align}
where the above inequality holds by the above inequalities~(\ref{eq:H8}) and~(\ref{eq:H11}). 

Let $0<\gamma\leq \frac{\lambda \mu^2}{108 L^2_f}$, we have $\frac{\gamma\eta}{6} - \frac{9\gamma^2\eta L^2_f}{\lambda\mu^2}\geq \frac{\gamma\eta}{12}$. Let $\gamma \leq \frac{\eta(1-\nu^2)}{12(H_t(1+\nu^2)+
2\theta_{x,t}\eta)}$, we have $H_t\frac{2(1+\nu^2)}{1-\nu^2}\gamma^2+\frac{4\theta_{x,t}\eta}{1-\nu^2}\gamma^2\leq \frac{\gamma\eta}{12}$.
Further let $\theta_{x,t}= \frac{2}{(1-\nu^2)\eta}\big( (3+\nu^2)H_t + 9\gamma\eta L_f^2+\frac{\gamma\eta L^2}{2}+\frac{288\gamma\eta L^4_f}{\mu^2} \big)$ and $\theta_{x,t}\leq \theta_{x,t-1}$ for all $t\geq1$,
we have $\theta_{x,t}(1-\frac{(1-\nu^2)\eta}{2})+(3+\nu^2)H_t \leq \theta_{x,t-1} - 9\gamma\eta L_f^2- \frac{\gamma\eta L^2}{2}-\frac{288\gamma\eta L^4_f}{\mu^2}$.
Meanwhile, let  $\theta_{y,t}= \frac{2}{(1-\nu^2)\eta}\big( (3+\nu^2)H_t + 18\gamma\eta L_f^2 + \frac{288\gamma\eta L^4_f}{\mu^2} \big)$ and $\theta_{y,t}\leq \theta_{y,t-1}$ for all $t\geq1$,
we have $\theta_{y,t}(1-\frac{(1-\nu^2)\eta}{2})+(3+\nu^2)H_t\leq \theta_{y,t-1} - 18\gamma\eta L_f^2 - \frac{288\gamma\eta L^4_f}{\mu^2}$.
According to the above inequality (\ref{eq:H18}), then we have
\begin{align} \label{eq:H19}
\Omega_{t+1} & \leq   \mathbb{E}_{t+1}\Big[ F(\bar{x}_t) + \frac{72\gamma L^2_f}{\lambda\mu^2}\big(F(\bar{x}_t) -f(\bar{x}_t,\bar{y}_t)\big) -\frac{\gamma\eta}{12m}\sum_{i=1}^m\|w^i_{x,t}-\nabla F(x^i_t)\|^2-\frac{\gamma\eta L^2}{4}\frac{1}{m}\sum_{i=1}^m\|x^i_{t}-\bar{x}_{t}\|^2   \nonumber \\
 & \qquad - \frac{\gamma\eta}{12} \frac{1}{m}\sum_{i=1}^m\|w^i_{x,t}\|^2 + \Delta_t +
 \frac{\gamma\eta}{12} \frac{1}{m}\sum_{i=1}^m\|w^i_{x,t-1}\|^2 \nonumber \\
& \qquad + (\rho_{x,t-1} - \frac{9\gamma\eta}{2})\mathbb{E} \|\bar{u}_{x,t-1} -\overline{\nabla_x f(x_{t-1},y_{t-1})}\|^2 + (\rho_{y,t-1}- \frac{144\gamma\eta L^2_f}{\mu^2})\mathbb{E} \|\bar{u}_{y,t-1} -\overline{\nabla_y f(x_{t-1},y_{t-1})}\|^2 \nonumber \\
& \qquad  + \varrho_{x,t-1}\frac{1}{m}\sum_{i=1}^m\mathbb{E} \|u^i_{x,t-1} - \nabla_x f^i(x^i_{t-1},y^i_{t-1})\|^2 + \varrho_{y,t-1}\frac{1}{m}\sum_{i=1}^m\mathbb{E} \|u^i_{y,t-1} - \nabla_y f^i(x^i_{t-1},y^i_{t-1})\|^2 \nonumber \\
& \qquad + \big( \theta_{x,t-1} - 9\gamma\eta L_f^2- \frac{\gamma\eta L^2}{2}-\frac{288\gamma\eta L^4_f}{\mu^2}\big)\frac{1}{m}\sum_{i=1}^m\|x^i_{t-1}-\bar{x}_{t-1}\|^2 \nonumber \\
& \qquad + \big(\theta_{y,t-1} - 18\gamma\eta L_f^2 - \frac{288\gamma\eta L^4_f}{\mu^2}\big)\frac{1}{m}\sum_{i=1}^m\|y^i_{t-1}-\bar{y}_{t-1}\|^2  +(\vartheta_{x,t-1}-\frac{3\gamma\eta}{4})\frac{1}{m}\sum_{i=1}^m\|w^i_{x,t-1}-\bar{w}_{x,t-1}\|^2 \nonumber \\
& \qquad + \big(\nu\vartheta_{y,t} + H_t\frac{4(1+\nu^2)}{1-\nu^2}\lambda^2 + \frac{8\theta_{y,t}\eta\lambda^2}{1-\nu^2} \big)\frac{1}{m}\sum_{i=1}^m\|w^i_{y,t-1}-\bar{w}_{y,t-1}\|^2 \nonumber \\
& \qquad + \big(H_t\frac{4(1+\nu^2)}{1-\nu^2}\lambda^2 + \frac{8\theta_{y,t}\eta\lambda^2}{1-\nu^2}\big)\|\bar{w}_{y,t-1}\|^2 \Big].
\end{align}

Let $\vartheta_{y,t}=\frac{4\lambda^2}{(1-\nu)(1-\nu^2)}\big(H_t(1+\nu^2) + 2\theta_{y,t}\eta\big)$ and
$\vartheta_{y,t} \leq \vartheta_{y,t-1}$ for all $t\geq 1$, we have $\nu\vartheta_{y,t} + H_t\frac{4(1+\nu^2)}{1-\nu^2}\lambda^2 + \frac{8\theta_{y,t}\eta\lambda^2}{1-\nu^2} = \vartheta_{y,t} -
(1-\nu)\vartheta_{y,t} + H_t\frac{4(1+\nu^2)}{1-\nu^2}\lambda^2 + \frac{8\theta_{y,t}\eta\lambda^2}{1-\nu^2} \leq \vartheta_{y,t-1}$. Let $\lambda \leq \frac{3L_f}{\mu}\sqrt{\frac{\gamma\eta(1-\nu^2)}{2H_t(1+\nu^2)+4\theta_{y,t}\eta}}$, we have $H_t\frac{4(1+\nu^2)}{1-\nu^2}\lambda^2 + \frac{8\theta_{y,t}\eta\lambda^2}{1-\nu^2} \leq \frac{18\gamma L^2_f\eta}{\mu^2}$. 
Then we have
\begin{align} \label{eq:H20}
\Omega_{t+1} & \leq   \mathbb{E}_{t+1}\Big[ F(\bar{x}_t) + \frac{72\gamma L^2_f}{\lambda\mu^2}\big(F(\bar{x}_t) -f(\bar{x}_t,\bar{y}_t)\big) -\frac{\gamma\eta}{12m}\sum_{i=1}^m\|w^i_{x,t}-\nabla F(x^i_t)\|^2-\frac{\gamma\eta L^2}{4}\frac{1}{m}\sum_{i=1}^m\|x^i_{t}-\bar{x}_{t}\|^2   \nonumber \\
 & \qquad - \frac{\gamma\eta}{12} \frac{1}{m}\sum_{i=1}^m\|w^i_{x,t}\|^2 + \Delta_t +
 \frac{\gamma\eta}{12} \frac{1}{m}\sum_{i=1}^m\|w^i_{x,t-1}\|^2 + \frac{18\gamma L^2_f\eta}{\mu^2}\|\bar{w}_{y,t-1}\|^2 \nonumber \\
& \qquad + (\rho_{x,t-1} - \frac{9\gamma\eta}{2})\mathbb{E} \|\bar{u}_{x,t-1} -\overline{\nabla_x f(x_{t-1},y_{t-1})}\|^2 + (\rho_{y,t-1}- \frac{144\gamma\eta L^2_f}{\mu^2})\mathbb{E} \|\bar{u}_{y,t-1} -\overline{\nabla_y f(x_{t-1},y_{t-1})}\|^2 \nonumber \\
& \qquad  + \varrho_{x,t-1}\frac{1}{m}\sum_{i=1}^m\mathbb{E} \|u^i_{x,t-1} - \nabla_x f^i(x^i_{t-1},y^i_{t-1})\|^2 + \varrho_{y,t-1}\frac{1}{m}\sum_{i=1}^m\mathbb{E} \|u^i_{y,t-1} - \nabla_y f^i(x^i_{t-1},y^i_{t-1})\|^2 \nonumber \\
& \qquad + \big( \theta_{x,t-1} - 9\gamma\eta L_f^2- \frac{\gamma\eta L^2}{2}-\frac{288\gamma\eta L^4_f}{\mu^2}\big)\frac{1}{m}\sum_{i=1}^m\|x^i_{t-1}-\bar{x}_{t-1}\|^2 \nonumber \\
& \qquad + \big(\theta_{y,t-1} - 18\gamma\eta L_f^2 - \frac{288\gamma\eta L^4_f}{\mu^2}\big)\frac{1}{m}\sum_{i=1}^m\|y^i_{t-1}-\bar{y}_{t-1}\|^2
 \nonumber \\
& \qquad +(\vartheta_{x,t-1}-\frac{3\gamma\eta}{4})\frac{1}{m}\sum_{i=1}^m\|w^i_{x,t-1}-\bar{w}_{x,t-1}\|^2 + \vartheta_{y,t-1} \frac{1}{m}\sum_{i=1}^m\|w^i_{y,t-1}-\bar{w}_{y,t-1}\|^2\Big] \nonumber \\
& = \Omega_t + \Delta_t - \frac{\gamma\eta}{12m}\sum_{i=1}^m\|w^i_{x,t}-\nabla F(x^i_t)\|^2-\frac{\gamma\eta L^2}{4}\frac{1}{m}\sum_{i=1}^m\|x^i_{t}-\bar{x}_{t}\|^2 - \frac{\gamma\eta}{12} \frac{1}{m}\sum_{i=1}^m\|w^i_{x,t}\|^2.
\end{align}

According to the above inequality~(\ref{eq:H20}), we can get  
\begin{align} \label{eq:H21}
 \frac{1}{m}\sum_{i=1}^m\|w^i_{x,t}-\nabla F(x^i_t)\|^2+\frac{ L^2}{m}\sum_{i=1}^m\|x^i_{t}-\bar{x}_{t}\|^2 + \frac{1}{m}\sum_{i=1}^m\|w^i_{x,t}\|^2 \leq \frac{12(\Omega_t - \Omega_{t+1})}{\gamma\eta} + \frac{12\Delta_t}{\gamma\eta}.
\end{align}

Let $\mathcal{M}^i_t = \|w^i_{x,t}-\nabla F(x^i_t)\| + \|w^i_{x,t}\| + L\|x^i_{t}-\bar{x}_{t}\| $, we have
\begin{align}
 \mathcal{M}^i_t & = \|w^i_{x,t}-\nabla F(x^i_t)\| + \|w^i_{x,t}\| + L\|x^i_{t}-\bar{x}_{t}\| \nonumber \\
 & \geq \|\nabla F(x^i_t)\|  + L\|x^i_{t}-\bar{x}_{t}\| \nonumber \\
 & \geq \|\nabla F(x^i_t)\|  + \|\nabla F(x^i_{t})-\nabla F(\bar{x}_{t})\| \geq \|\nabla F(\bar{x}_{t})\|
\end{align}
According to the above inequality (\ref{eq:H21}), then we can obtain
\begin{align} \label{eq:H22}
\frac{1}{T}\sum_{t=1}^T\frac{1}{m}\sum_{i=1}^m\mathbb{E}[\mathcal{M}^i_t]^2  &
\leq \frac{1}{T}\sum_{t=1}^T\big[\frac{3}{m}\sum_{i=1}^m\|w^i_{x,t}-\nabla F(x^i_t)\|^2 +\frac{3L^2}{m}\sum_{i=1}^m\|x^i_t-\bar{x}_t\|^2 + \frac{3}{m}\sum_{i=1}^m\|w^i_{x,t}\|^2 \big] \nonumber \\
& \leq \frac{1}{T}\sum_{t=1}^T\frac{36(\Omega_t-\Omega_{t+1})}{\gamma\eta} \leq \frac{36(\Omega_1-F^*)}{T\gamma\eta} + \frac{36\sum_{t=1}^T\Delta_t}{T\gamma\eta}.
\end{align}
By using the Cauchy-Schwarz inequality, we can further get 
 \begin{align} \label{eq:H23}
  \frac{1}{T}\sum_{t=1}^T\frac{1}{m}\sum_{i=1}^m\mathbb{E}\big[\mathcal{M}^i_t \big] \leq \sqrt{\frac{1}{T}\sum_{t=1}^T\frac{1}{m}\sum_{i=1}^m\mathbb{E}[\mathcal{M}^i_t]^2} \leq \frac{1}{T}\sum_{t=1}^T\frac{36(\Omega_t-\Omega_{t+1})}{\gamma\eta} \leq \frac{6\sqrt{\Omega_1-F^*}}{\sqrt{T\gamma\eta}} + \frac{6\sqrt{\sum_{t=1}^T\Delta_t}}{\sqrt{T\gamma\eta}}.
 \end{align}
Thus we can obtain 
\begin{align}
  \frac{1}{T}\sum_{t=1}^T \|\nabla F(\bar{x}_t)\|
  \leq  \frac{1}{T}\sum_{t=1}^T\frac{1}{m}\sum_{i=1}^m\mathbb{E}\big[\mathcal{M}^i_t \big]  \leq \frac{6\sqrt{\Omega_1-F^*}}{\sqrt{T\gamma\eta}} + \frac{6\sqrt{\sum_{t=1}^T\Delta_t}}{\sqrt{T\gamma\eta}}.
\end{align}

Let $\alpha_t =\beta_t = O(\frac{1}{T^{2/3}})$ and $\eta_t=\eta=O(\frac{1}{T^{1/3}})$ for all $t\geq1$, and
$\gamma=O(1)$ and $\lambda=O(1)$, we have
\begin{align}
 & \rho_{x,t} = \frac{9\gamma\eta}{2\alpha_t} = O(T^{1/3}), \quad \rho_{y,t} = \frac{144\gamma\eta L^2_f}{\mu^2\beta_t}=O(T^{1/3}), \quad \vartheta_{x,t}=\frac{\gamma\eta}{1-\nu}=O(T^{-1/3}), \quad \vartheta_{y,t}=O(T^{-1/3}) \nonumber \\
 & \varrho_{x,t} =\frac{4\vartheta_{x,t}\nu^2}{1-\nu}=O(T^{-1/3}), \quad \varrho_{y,t} = O(T^{-1/3}),
\end{align}
and 
\begin{align}
 \Delta_t = \frac{2\rho_{x,t}\alpha_{t}^2\sigma^2}{m}+ \frac{2\rho_{y,t}\beta_{t}^2\sigma^2}{m}+ 2\varrho_{x,t}\alpha_{t}^2\sigma^2 + 2\varrho_{y,t}\beta_{t}^2\sigma^2 + \frac{4\nu^2\vartheta_{x,t} \alpha^2_{t}\sigma^2}{1-\nu} + \frac{4 \nu^2\vartheta_{y,t}\beta^2_{t}\sigma^2}{1-\nu}=O(\frac{\sigma^2}{T}).
\end{align}
Further let $\rho_{x,0} = \frac{9\gamma\eta}{2}$, $\rho_{y,0} = \frac{144\gamma\eta L^2_f}{\mu^2}$,
$\varrho_{x,0}=\varrho_{y,0}=0$, $x_0^1=\tilde{x}_0^1=\cdots=x_0^m=\tilde{x}_0^m$, $y_0^1=\tilde{y}_0^1=\cdots=y_0^m=\tilde{y}_0^m$, $w^1_{x,0}=\cdots=w^m_{x,0}=0$ and $w^1_{y,0}=\cdots=w^m_{y,0}=0$,
 we have
\begin{align} 
 \Omega_1 = F(\bar{x}_1) + \frac{72\gamma L^2_f}{\lambda\mu^2}\big(F(\bar{x}_1) 
 -f(\bar{x}_1,\bar{y}_1)\big) =O(1).
\end{align}
Then we can get 
\begin{align}
  \frac{1}{T}\sum_{t=1}^T \|\nabla F(\bar{x}_t)\|
  \leq O(\frac{1}{T^{1/3}}+\frac{\sigma^2}{T^{1/3}}) \leq \epsilon.
\end{align}
Without loss of generality, let $T= O(\epsilon^{-3})$, then we can obtain the gradient (SFO) complexity $1\cdot T=O(\epsilon^{-3})$ of our DM-GDA algorithm.

\end{proof}

\end{document}